\documentclass[article]{siamart171218}
\usepackage{pgfplots,subfig}
\usepackage{tikz}
\usepackage{ifthen}
\pgfplotsset{compat=1.5}
\usepackage{graphicx}
\usepackage{caption}
\usepackage{float}
\usepackage{amsopn,amssymb,amsmath}
\usepackage{mathrsfs}
\usepackage{algorithmic}
\usepackage{array,booktabs,tabularx,makecell, blkarray}
\usepackage[titletoc,toc,title]{appendix}
\pgfplotsset{every tick label/.append style={font=\footnotesize}}

\usepackage{etoolbox}
\providetoggle{comments}\settoggle{comments}{true}
\iftoggle{comments}{
    \usepackage{todonotes}
    \newcommand\ed[1]{\todo[color=red!10,inline]{ED:~#1}}
    \newcommand\lc[1]{\todo[color=orange!10,inline]{LC:~#1}}
    \newcommand\zx[1]{\todo[color=green!10,inline]{ZX:~#1}}
    \newcommand\pe[1]{\todo[color=blue!10,inline]{PE:~#1}}
    \newcommand\fh[1]{\todo[color=brown!10,inline]{FHR:~#1}}
    \newcommand\yh[1]{\todo[color=cyan!10,inline]{YH:~#1}}
    \newcommand\ca[1]{\todo[color=magenta!10,inline]{CA:~#1}}
}{
    \newcommand\ed[1]{}
    \newcommand\lc[1]{}
    \newcommand\zx[1]{}
    \newcommand\pe[1]{}
    \newcommand\fh[1]{}
    \newcommand\yh[1]{}
    \newcommand\ca[1]{}
}

\makeatletter
\newcommand{\@algocf@capt@plain}{above}
\makeatother

\graphicspath{{Images/}}

\DeclareMathOperator{\dist}{dist}
\DeclareMathOperator{\dr}{dr}
\newcommand\nobreakpar{\par\nobreak\@afterheading}

\newcommand{\wtX}{\widetilde{X}}
\newcommand{\wtY}{\widetilde{Y}}
\newcommand{\whK}{\widehat{K}}
\newcommand{\whX}{\widehat{X}}
\newcommand{\whY}{\widehat{Y}}
\newcommand{\Xinit}{X_{init}}
\newcommand{\Yinit}{Y_{init}}
\renewcommand{\Xinit}{X^\circ}
\renewcommand{\Yinit}{Y^\circ}
\newcommand{\Xhat}{\widehat{X}}
\newcommand{\Yhat}{\widehat{Y}}
\newcommand{\whKxy}{\whK_{\Xinit,\Yinit}}
\newcommand{\Wx}{W_{\Xinit,\Xinit}}
\newcommand{\Wy}{W_{\Yinit,\Yinit}}

\newcommand{\Kfun}{\mathscr{K}}
\newcommand{\Xdom}{\bf{X}}
\newcommand{\Ydom}{\bf{Y}}

\newcommand{\si}{{\text{SI}}}

\newcommand{\KXY}{K_{X,Y}}
\newcommand{\KXhatY}{K_{\Xhat,Y}}
\newcommand{\KXYhat}{K_{X,\Yhat}}
\newcommand{\KXhatYhat}{K_{\Xhat,\Yhat}}
\newcommand{\Kinit}{K_{\Xinit,\Yinit}}

\usepackage{hyperref}

\def\equationautorefname~#1\null{Eq.~(#1)\null}

\makeatletter
\newcommand\Autoref[1]{\@first@ref#1,@}
\def\@throw@dot#1.#2@{#1}
\def\@set@refname#1{
    \edef\@tmp{\getrefbykeydefault{#1}{anchor}{}}%
    \xdef\@tmp{\expandafter\@throw@dot\@tmp.@}%
    \ltx@IfUndefined{\@tmp autorefnameplural}%
         {\def\@refname{\@nameuse{\@tmp autorefname}s}}%
         {\def\@refname{\@nameuse{\@tmp autorefnameplural}}}%
}
\def\@first@ref#1,#2{%
  \ifx#2@\autoref{#1}\let\@nextref\@gobble
  \else%
    \@set@refname{#1}
    \@refname~\ref{#1}
    \let\@nextref\@next@ref
  \fi%
  \@nextref#2%
}
\def\@next@ref#1,#2{%
   \ifx#2@ and~\ref{#1}\let\@nextref\@gobble
   \else, \ref{#1}
   \fi%
   \@nextref#2%
}
\makeatother

\title{Low-rank kernel matrix approximation \\
       using Skeletonized Interpolation \\
       with Endo- or Exo-Vertices}

\author{
   Zixi Xu\thanks{Department of Mechanical Engineering, 
                   Stanford University (\email{zixixu@stanford.edu}, 
                   \email{darve@stanford.edu})}
    \and L{\'e}opold Cambier\thanks{Institute for Computational 
                 \& Mathematical Engineering, Stanford University 
                 (\email{lcambier@stanford.edu})}
    \and Francois-Henry Rouet\thanks{Livermove Software Technology
                 Corporation 
                 \email{fhrouet@lstc.com}, 
                 (\email{pierre@lstc.com}, 
                 \email{huang@lstc.com}, 
                 \email{cleve@lstc.com})}
    \and Pierre L'Eplattenier\footnotemark[3]
    \and Yun Huang\footnotemark[3]
    \and Cleve Ashcraft\footnotemark[3]
    \and Eric Darve\footnotemark[1]
}

\begin{document}

\maketitle

\begin{abstract}
The efficient compression of kernel matrices, for instance the
off-diagonal blocks of discretized integral equations, is a crucial
step in many algorithms.
In this paper, we study the application of Skeletonized
Interpolation to construct such factorizations.
In particular, we study four different strategies for selecting the
initial candidate pivots of the algorithm: Chebyshev grids, points
on a sphere, maximally-dispersed and random vertices.
Among them, the first two introduce new interpolation points
(exo-vertices) while the last two are subsets of the given clusters
(endo-vertices). We perform experiments using three real-world
problems coming from the multiphysics code LS-DYNA. The pivot
selection strategies are compared in term of quality (final rank)
and efficiency (size of the initial grid).
These benchmarks demonstrate that overall, maximally-dispersed
vertices provide an accurate and efficient sets of pivots for most
applications. It allows to reach near-optimal ranks while starting
with relatively small sets of vertices, compared to other strategies.
\end{abstract}


\section{Introduction} \label{Introduction}


Our goal is to compute low-rank approximations of matrices that come
from certain kernel functions
\begin{equation}
\label{eqn:}
\Kfun : \Xdom \times \Ydom
\end{equation}
where $\Xdom$ and $\Ydom$ are 2-D or 3-D regions of space.
From $\Xdom$ and $\Ydom$ we choose sets of points 
\begin{equation*}
X = \{x_1,\hdots,x_m\}
\qquad \text{and} \qquad
Y = \{y_1,\hdots,y_n\}
\end{equation*}
on which to evaluate the kernel. This yields a $m \times n$ matrix,
$\KXY$, where $K_{i,j} = \Kfun(x_i,y_j)$.
%
When $\Kfun$ is smooth over the domain surrounding $\Xdom$ and
$\Ydom$,  $\KXY$ typically has rapidly decaying singular values and can be well approximated by a low-rank matrix.

Matrices of this sort arise naturally in many applications. 
A typical electromagnetic application may have at a boundary
integral equation 
\begin{equation} \label{eq:EM}
\vec{\Phi}(x) 
= \frac{\mu_0}{4\pi} 
  \int_{\Gamma_{\Ydom}} \Kfun(x,y)\, \vec{k}(y)\,dy
= \frac{\mu_0}{4\pi} 
  \int_{\Gamma_{\Ydom}}\frac{1}{\|x-y\|_2}\, \vec{k}(y)\,dy
  \qquad \text{for}\ x \in \Gamma_{\Xdom}
\end{equation}
where $\Gamma_{\Xdom}$ and $\Gamma_{\Ydom}$ are the boundaries of
$\Xdom$ and $\Ydom$. 
Here $\vec{\Phi}$ is the
magnetic vector potential, $\mu_0$ is the electrical permeability and
$\vec{k}$ is the surface current introduced as an intermediate variable \cite{EM_BEM2009}.


In acoustics, a boundary integral equation based on Helmholtz equations
is complex
\begin{equation} \label{eq:Acoustics}
\frac{1}{2} p(x)= -\int_{\Gamma_y} \left( i \rho \omega v_n(y) \mathscr{G} + p(y) \frac{\partial \mathscr{G}}{\partial n} \right)dy.
\end{equation}
where $p(x)$, $p(y)$ are acoustic pressures, 
$v_n(y)$ is the normal velocity,
$\rho$ is the fluid density, $\omega$ is the round frequency 
\cite{wrobel2002boundary}. 

In equations~(\ref{eq:EM}) and~(\ref{eq:Acoustics}), the kernel
matrices are given by
\begin{equation*}
\Kfun(x,y)       = \frac{1}{r}          
\qquad \text{and} \qquad
\mathscr{G}(x,y) = \frac{1}{r} e^{-ikr}
\end{equation*}
where $r = \|x-y\|_2$.
Other BEM kernel functions of interest include $\ln(r)$ and $1/r^2$.


The discretized forms of 
\autoref{eq:EM} and
\autoref{eq:Acoustics} 
are linear systems
\begin{equation} \label{eq:discretized}
    K_{N,N} \, u_{N} = b_{N}
\end{equation}
where the $|N| \times |N|$ matrix $K_{N,N}$ is dense, 
with solution $u_N$ and right hand side $b_N$.
A submatrix $K_{X,Y}$, where $X \subset N$ and $Y \subset N$,
will be dense, but if the subsets $X$ and $Y$ are chosen well,
$K_{X,Y}$ will have small numerical rank.
\par
In this case, we have 
$K_{X,Y} \approx U_{X,\alpha} V_{Y,\alpha}^{\top}$
where typically $|X|$ and $|Y|$ are $O(100)$ 
and $|\alpha|$ is $O(1)$ or $O(10)$.
Low rank storage is less than dense storage, 
$(|X|+|Y|)|\alpha|$ vs $|X|\,|Y|$.
A matrix-vector multiply with $K_{X,Y}$, 
used in an iterative solution of \autoref{eq:discretized}, 
takes fewer operations than if $K_{X,Y}$ were a dense matrix.
\par
One approach is to form the dense $K_{X,Y}$
and then compute a rank-revealing QR factorization.
If $|X| = |Y| = n$ and $|\alpha| = r$, this takes $n^2$ kernel
evaluations to form $K_{X,Y}$ and $O(n^2r)$ operations
for the RRQR.
Our goal is to reduce these two costs, to $O(nr)$ kernel evaluations
and $O(nr^2)$ or even $O(nr)$ linear algebra operations.


A new low-rank approximation approach for the kernel matrix is
proposed in this paper. For illustration purpose, this paper takes
the kernel functions $1 / r$ and $1 / r^2$ ($r=\|x-y\|_2$) as
examples to show how the approach works. Note that although our
approach only focuses on each single sub-matrix extracted from the
fully dense matrix, it can be integrated as the building block of a
more complex multi-level method using hierarchical matrices for
instance.

\ifthenelse{1=1}{}{
where $X$ here denotes all the vertices in the mesh.
\fh{This look weird to me. It can't be a single pair (X,Y). Is there a sum over Y's?}
\lc{Is this better ?}
However, $K_{X,X}$ in \autoref{eq:discretized} as well as the resultant
\lc{Updated}
kernel matrices from \autoref{tab:BEMKernel}
are fully dense such that the time and storage of computing the kernel matrix will grow as
$O(N^2)$ where $N$ is the size of the matrix.
Although these kernel matrices are usually
not low-rank, one can find subsets of points $X$ and $Y$ corresponding to well-separated domains
such that $\KXY$ is low-rank.
Performing an efficient low-rank approximation on
$\KXY$ can lead to significant storage savings. In particular, it would be desirable to be able
to factor $\KXY$ into a low-rank representation $\KXY \approx UV^\top$ \emph{without}
having to form the matrix.

A new low-rank approximation approach for the kernel matrix is
proposed in this paper. For illustration purpose, this paper takes
the kernel functions $1 / r$ and $1 / r^2$ ($r=\|x-y\|_2$) as
examples to show how the approach works. Note that although our
approach only focuses on each single sub-matrix extracted from the
fully dense matrix, it can be integrated as the building block of a
more complex multi-level method using hierarchical matrices for
instance.

\begin{table}[!ht]
\caption{Typical BEM problems and their corresponding kernel functions}
\label{tab:BEMKernel}
\begin{center}
\begin{tabular}{ll}
\toprule
\textbf{BEM problems}  &   \textbf{Kernel functions ($r=\|x-y\|_2$)}                      \\
\midrule
Electromagnetics       &   $1/r$                          \\
\addlinespace[0.7em]
Acoustics              &   $\exp(-ikr) / r$           \\
\addlinespace[0.7em]
2D potential           &   $1/r$, $\ln r$  \\
\addlinespace[0.7em]
Elastostatics          &   $1/r$, $1 / r^2$ \\
\bottomrule
\end{tabular}
\end{center}
\end{table}
} 

For convenience, \autoref{tab:TableOfNotation} summarizes the
notations we use throughout the paper. \autoref{tab:Acronyms}
summarizes the terminology (acronyms and algorithms).

\begin{table}[!ht]\caption{Notation}
\label{tab:TableOfNotation}
\centering
\begin{tabular}{l p{260pt}}
\toprule
$\Kfun$ & The kernel function\\
$\Xdom$, $\Ydom$ & Two regions in space \\
$X$, $Y$ & Sets of discretization points in subdomains of a geometry\\
$\Xinit$, $\Yinit$ & Initial interpolation points based on $X$ and $Y$\\
$\Xhat$, $\Yhat$ & Subsets of $X$ and $Y$ chosen by Skeletonized
   Interpolation \\
      & used to build the low-rank approximation\\
$\KXY$ & The kernel matrix defined as $(\KXY)_{ij} = \Kfun(x_i,y_j)$ \\
$\Wx$, $\Wy$ & Corresponding weight matrices for $\Xinit$ and $\Yinit$\\
$r_0$ & Size of $\Xinit$ or $\Yinit$ \\
$r_1$ & Size of $\Xhat$ or $\Yhat$\\
$\text{dr}(i,j)$ & The distance ratio between two pairs $i$ and $j$ of clusters\\
$\varepsilon$ & The tolerance used in the RRQR algorithms\\
$\varepsilon^*$ & Relative Frobenius-norm error used in the experiments\\
\bottomrule
\end{tabular}
\end{table}

\begin{table}[!ht]\caption{Terminology}
\label{tab:Acronyms}
\centering
\begin{tabular}{l p{280pt}}
\toprule
SI & Skeletonized Interpolation \\
MDV & Maximally-dispersed vertices \\
RRQR & Rank-revealing QR\\
{\sf SI-Chebyshev} & SI algorithm using Chebyshev nodes for $\Xinit$ and $\Yinit$\\
{\sf SI-MDV} & SI algorithm using maximally-dispersed vertices for $\Xinit$, $\Yinit$\\
{\sf SI-sphere} & SI algorithm using points on a bounding sphere for $\Xinit$, $\Yinit$\\
{\sf SI-random} & SI algorithm using random vertices for $\Xinit$ and $\Yinit$\\
${\sf GenInitSet}$ & An algorithm that returns, for a set of vertices $X$, a set $\Xinit$ of a given maximum size and the corresponding weights $\Wx$\\
\bottomrule
\end{tabular}
\end{table}

\subsection{Previous work}

The efficient solution of boundary integral equations [e.g.,
Eqns.~\eqref{eq:EM} and \eqref{eq:Acoustics}] has
been extensively studied. 
Multiple algorithms have been proposed to efficiently compute 
low-rank approximations for the off-diagonal blocks $\KXY$ of 
the full kernel matrix.

A variety of analytical expansion methods are
based on the expansion properties of the kernel, 
The Fast Multipole Method
(FMM) \cite{Nishimura2002, darve2000fast1, darve2000fast2} is based
on a series expansion of the fundamental solution. 
It was first proposed by \cite{Rokhlin1985} for the Laplacian, 
and then a diagonal version for the Helmholtz operator 
\cite{Greengard1997}. 
FMM can accelerate matrix-vector products which are coupled 
with an iterative method.
More recently, \cite{Fong2009TheBF} proposed a kernel-independent
FMM based on the interpolation of the kernel.
\par
In \cite{chew2004boundaryEM}, Chew et.\ al.\ 
presented their benchmark results and scaling studies
for the multilevel fast multipole algorithm (MLFMA) 
and the fast inhomogeneous plane wave algorithm (FIPWA). 
A hybrid version of MLFMA only requires a fraction of CPU time 
and memory of the traditional MLFMA.
FIPWA largely reduces the CPU time thus making simulations 
feasible with over one million particles.

In addition to the FMM, there exists other techniques that
approximate the kernel function to explicitly compute
the low-rank factorization of the kernel submatrices. 
The Panel Clustering method \cite{Hackbusch1989} provides a kernel
function approximation using Taylor Series. 
\cite{Yarvin1998} (and similarly \cite{Borm2004} and \cite{wu2014}
      in Fourier space) takes advantage of the interpolation of
$\Kfun(x,y)$ over $X \times Y$ to build a low-rank expansion $\KXY =
S_{X,\widetilde{X}}K_{\widetilde{X},\widetilde{Y}}
T_{Y,\widetilde{Y}}^\top$ and uses the SVD to recompress that
expansion further.
Also, Barnes and Hut \cite{barnes1986Nature} compute mutual forces
in particle systems using a center of mass approximation with a
special far-field separation or ``admissibility'' criterion. 
This is similar to an order-1 multipole expansion (where one matches
multipole terms of order up to 1).

However, these analytical methods are often limited to some specific
types of kernel functions, or have complexities larger than $O(nr)$.
Adaptive Cross Approximation (ACA) \cite{Bebendorf2003, Bebendorf2000} 
computes individual rows and columns of the matrix.
However, it provides few accuracy guarantees and its termination
criterion is inaccurate in some cases. Furthermore, the method is
not very efficient because it proceeds column by column and row
by row instead of using matrix-matrix operations (level 3 BLAS
subprograms) that have a more efficient memory access
pattern.

Finally, Bebendorf \cite{Bebendorf2000} proposes the form
\begin{equation} \label{eq:Bebendorf}
\KXY = K_{X,\widetilde{Y}} K_{\widetilde{X},\widetilde{Y}}^{-1} K_{\widetilde{X},Y}
\end{equation}
where $\widetilde{X}$ and $\widetilde{Y}$ are interpolation points
built iteratively from $\Kfun(x,y)$. Our method explores this
interpolation idea; however, we choose the interpolation nodes in a
very different way. While \cite{Bebendorf2000} builds the nodes one
by one in an adaptive fashion, we start from a predefined grid of
nodes and then sample this grid using a RRQR factorization. This
leads to a fast and very easy algorithm.




\subsection{Skeletonized Interpolation} \label{SI-intro}

We now introduce our approach.
Consider the kernel matrix $\KXY$ having an exact rank $r$. 
Let $\wtX \subseteq X$ and $\wtY \subseteq Y$ be sampling points
such that $|\wtX| = |\wtY| = r$ and $K_{\wtX,\wtY}$ is nonsingular.
Write the kernel matrix in block form and factor as
\begin{align}
\notag
K_{X,Y} 
&= 
\begin{bmatrix}
K_{\wtX,\wtY}           & K_{\wtX,Y\setminus\wtY}           \\
K_{X\setminus\wtX,\wtY} & K_{X\setminus\wtX,Y\setminus\wtY} \\
\end{bmatrix}
\\ \label{eq:endo}
&= 
\begin{bmatrix}
I_{\wtX,\wtY}           \\
L_{X\setminus\wtX,\wtY} \\
\end{bmatrix}
K_{\wtX,\wtY}
\begin{bmatrix}
I_{\wtX,\wtY} & U_{\wtX,Y\setminus\wtY} \\
\end{bmatrix}
+
\begin{bmatrix}
0_{\wtX,\wtY}           & 0_{\wtX,Y\setminus\wtY}           \\
0_{X\setminus\wtX,\wtY} & E_{X\setminus\wtX,Y\setminus\wtY} \\
\end{bmatrix}
\end{align}
where
$L_{X\setminus\wtX,\wtY} = K_{X\setminus\wtX,\wtY} K_{\wtX,\wtY}^{-1}$
and
$U_{\wtX,Y\setminus\wtY} = K_{\wtX,\wtY}^{-1} K_{\wtX,Y\setminus\wtY}$.
If $K_{X,Y}$ has numerical rank $r$, as does $K_{\wtX,\wtY}$, then
then error matrix $E_{X\setminus\wtX,Y\setminus\wtY}$ is zero,
and we have this low rank factorization of $K_{X,Y}$.
\begin{align}
\label{eq:endo-schur}
K_{X,Y}
&=
\begin{bmatrix}
I_{\wtX,\wtX}           \\
L_{X\setminus\wtX,\wtX} \\
\end{bmatrix}
K_{\wtX,\wtY}
\begin{bmatrix}
I_{\wtY,\wtY} & U_{\wtY,Y\setminus\wtY} \\
\end{bmatrix}
= K_{X,\wtY} K_{\wtX,\wtY}^{-1} K_{\wtX,Y} 
\end{align}
We call sampled degrees of freedom $\wtX \subseteq X$ 
and $\wtY\subseteq Y$ {\bf endo-vertices}, for they are internal
to the domains $X$ and $Y$.
\par
In contrast, consider {\bf exo-vertices}, where
$\wtX$ and $\wtY$ are chosen from outside the sets $X$ and $Y$,
$K_{\wtX,\wtY}$ is square and nonsingular.
Form the large kernel matrix and factor.
\begin{align}
\label{eq:exo}
K_{X\cup\wtX,Y\cup\wtY}
&= 
\begin{bmatrix}
K_{\wtX,\wtY} & K_{\wtX,Y} \\
K_{X,\wtY}    & K_{X,Y}    \\
\end{bmatrix}
= 
\begin{bmatrix}
I_{\wtX,\wtX} \\
L_{X,\wtX} \\
\end{bmatrix}
K_{\wtX,\wtY} 
\begin{bmatrix}
I_{\wtY,\wtY} & U_{Y,\wtY} \\
\end{bmatrix}
+
\begin{bmatrix}
0 & 0 \\
0 & E_{X,Y} \\
\end{bmatrix}
\end{align}
If the error matrix $E_{X,Y}$ is small in norm, then we have this
low rank representation for $K_{X,Y}$. 
\begin{equation}
\label{eq:exo-schur}
K_{X,Y} 
= L_{X,\wtX} K_{\wtX,\wtY} U_{\wtY,Y}
= K_{X,\wtY} K_{\wtX,\wtY}^{-1} K_{\wtX,Y}
\end{equation}
where
$L_{X,\wtX} = K_{X,\wtY} K_{\wtX,\wtY}^{-1}$
and
$U_{\wtY,Y} = K_{\wtX,\wtY}^{-1} K_{\wtX,Y}$.
\par
In both equations~\autoref{eq:endo-schur} and~\autoref{eq:exo-schur},
the large kernel matrix $K_{X,Y}$ as well as the smaller
$K_{\wtX,\wtY}$ matrices have log-linear singular values,
and so these matrices are ill-conditioned.
However, in \cite{LeopoldSI2017}, it was shown that if backward-stable 
algorithms are used to factor $K_{\wtX,\wtY}$, the product of the 
three matrices on the right-hand-sides can be computed accurately.
\par
There is a key difference between choosing endo-vertices and
exo-vertices.
The choice of endo-vertices $\wtX$ and $\wtY$ means that
for endo-,
$K_{\wtX,\wtY}$ must be able to capture the log-linear singular
values of $K_{X,Y}$.
For exo-,
$K_{\wtX,\wtY}$ must be able to capture the log-linear singular
values of the {\bf larger} kernel matrix $K_{\wtX\cup X,\wtY\cup Y}$.
The numerical rank of $K_{\wtX\cup X,\wtY\cup Y}$ will generally
be larger than the numerical rank of $K_{\wtX,\wtY}$.
We will see this to be the case in the experiments to follow.
\ifthenelse{1=1}{}{
Denote $\widetilde{X}$ a set of $r$ points sampled from
$X$ and $\widetilde{Y}$ from $Y$ such that
$K_{\widetilde{X},\widetilde{Y}}$ has rank $r$ and is invertible. In the following, we denote by $[\widetilde X, X]$ the list of vertices $X$ with $\widetilde X$ appended at the beginning. We consider the Schur complement
of \[K_{[\widetilde{X},X],[\widetilde{Y},Y]}\] when eliminating $\widetilde{X} \times \widetilde{Y}$. It is equal to
\[
    \KXY-K_{X,\widetilde{Y}}K_{\widetilde{X},\widetilde{Y}}^{-1}K_{\widetilde{X},Y}
\]
which is identically zero, since the rank is revealed. Or equivalently,
\begin{equation} \label{eq:schur}
    \KXY=K_{X,\widetilde{Y}}K_{\widetilde{X},\widetilde{Y}}^{-1}K_{\widetilde{X},Y}
\end{equation}

An extension of this idea is to select points $\widetilde{X}$ and $\widetilde{Y}$ outside the sets $X$ and $Y$. If the kernel \emph{function} $\Kfun(x, y)$ has an exact rank $r$, the above analysis stays valid.
However, the kernel function $\Kfun(x,y)$ and the kernel matrix $\KXY$ often have different ranks. In this case, the rank of $K_{[\widetilde X, X], [\widetilde Y, Y]}$ may be larger.


The problem now becomes how to choose sets $\widetilde{X}$ and 
$\widetilde{Y}$ in an easy, fast and accurate way. 
If we are able to do so, the low rank
approximation is then given by \autoref{eq:schur}.

In \autoref{eq:schur}, the central matrix
$K_{\widetilde{X},\widetilde{Y}}$ is ill-conditioned and its inverse
is inaccurate. 
However, in \cite{LeopoldSI2017}, it was shown that
if backward-stable algorithms are used, the product of the three
matrices on the right-hand-side of \autoref{eq:schur} can be
computed accurately.

In practice, if one has such a factorization, a (partially-pivoted)
LU factorization $\widetilde L \widetilde U \widetilde P =
K_{\widetilde{X}, \widetilde{Y}}$ is computed. Then, one stores the
following factorization
\[ \KXY = U V \quad \text{with} \quad U = K_{X, \widetilde Y}
\widetilde P^{-1} \widetilde U^{-1} \quad \text{and}\quad V =
   \widetilde L^{-1} K_{\widetilde{X}, Y} \]
} 

The problem now becomes how to choose sets $\widetilde{X}$ and 
$\widetilde{Y}$ in an easy, fast and accurate way. 
To select $\widetilde{X}$ and $\widetilde{Y}$, 
we start from larger sets of initial points $\Xinit \subseteq X$
and $\Yinit \subseteq Y$. 
We introduce different ways of choosing those points
in the next sections; we keep their exact definition unspecified at
the moment. 
\par
Using the method proposed in \cite{LeopoldSI2017}, we
perform two RRQR factorizations over 
\begin{equation*}
\whKxy = \Wx^{1/2} \Kinit\Wy^{1/2} 
\end{equation*}
and its transpose in order to choose an optimal set of rows and
columns, $\Xhat \subset \Xinit$ and $\Yhat \subset \Yinit$, for a
given tolerance $\varepsilon$. 
The diagonal matrices $\Wx$ and $\Wy$ contain
integration weights related to the choice of $\Xinit$ and $\Yinit$ (more details about this in \autoref{subsec:Discussion}).
Once $\Xhat$ and $\Yhat$ are selected, the resulting approximation
is given by 
\begin{equation}
\label{eqn:lowrankKXY}
\KXY \approx \KXYhat \KXhatYhat^{-1} \KXhatY 
\end{equation}
as indicated above.
\par
Let us consider the expense to create this low rank representation.
For simplicity, let $|X| = |Y| = n$,
$|\Xinit| = |\Yinit| = r_0$,
and $|\whX| = |\whY| = r_1$. 
\begin{itemize}
\item
There are $r_0^2$, $r_1^2$, $r_1n$ and $r_1n$ 
kernel function evaluations to compute
$K_{\Xinit,\Yinit}$, $K_{\whX,\whY}$, $K_{X\setminus\whX,\whY}$
and
$K_{\whX,Y\setminus\whY}$, respectively, for a total of 
$2r_1n + r_0^2 + r_1^2$.
\item
Two RRQR factorizations on $K_{\Xinit,\Yinit}$
for $8r_0^2 r_1$ linear algebra operations.
To invert $K_{\whX,\whY}$ requires $\frac{2}{3}r_1^3$ operations.
for a total of $8r_0^2 r_1 + \frac{2}{3}r_1^3$ operations.
\end{itemize}
\par
Ideally,
$\Xinit$ and $\Yinit$ are as small as possible, and $\Xhat$
and $\Yhat$ have (almost) the same sizes as the
$\varepsilon$-rank of $\KXY$.  The smaller those sets, the
less expensive the factorization is to compute, and the better
the sets.  In addition, as we discuss in \autoref{subsec:Discussion},
one can always further recompress a given low-rank
approximation. However, in this paper, we are interested in
the cost of the construction of the initial low-rank basis
based on the initial $\Xinit$ and $\Yinit$ and reduced
$\Xhat$ and $\Yhat$.


We call this method Skeletonized Interpolation. A high-level
description of the Skeletonized Interpolation algorithm is given in
\autoref{alg:genericSI}; ${\sf GenInitSet}$ denotes a given function
(algorithm) used to produce  $\Xinit$ and $\Yinit$; given a set $X$
or $Y$ and a given size $r_0$, it returns an initial set $\Xinit$ or
$\Yinit$ of size at most $r_0$, and the corresponding integration
weights.
\par
In this algorithm, as a rule of thumb, the tolerance $\varepsilon$
should be slightly smaller than the desired final accuracy. This
comes from the fact that the RRQR's are sub-optimal compared to the
SVD (in the sense that they typically lead to slightly larger ranks
for a given accuracy). 

\begin{algorithm}
\caption{Skeletonized Interpolation: $[\whX,\whY] = {\sf SI}(\Kfun, 
         X, Y, {\sf GenInitSet}, r_0, \varepsilon)$}
\label{alg:genericSI}
\begin{algorithmic}[1]
\REQUIRE{Kernel $\Kfun$, 
   discretization points $X$ and $Y$, 
   algorithm ${\sf GenInitSet}$, rank $r_0$,
   tolerance $\varepsilon$}
\STATE Calculate $(\Xinit, \Wx) = {\sf GenInitSet}(X, r_0)$ 
\STATE Calculate $(\Yinit, \Wy) = {\sf GenInitSet}(Y, r_0)$ 
\STATE Build temporary matrix \ 
   $T_{\Xinit,\Yinit} = 
   W_{\Xinit,\Xinit} K_{\Xinit,\Yinit} W_{\Yinit,\Yinit} $
   \label{SI_line2}
\STATE Perform truncated RRQR \ 
   $T_{\Xinit,\Yinit} P_{\Yinit,\Yinit} = 
    Q_{\Xinit,\alpha} R_{\alpha,\Yinit} + E_{\Xinit,\Yinit}$ \\
    \ \ \quad where $\|E_{\Xinit,\Yinit}\|_F \le \varepsilon$
    \label{SI_line3}
\STATE Perform truncated RRQR \ 
   $T_{\Xinit,\Yinit}^\top P_{\Xinit,\Xinit} = 
    Q_{\Yinit,\beta} R_{\beta,\Xinit} + E_{\Yinit,\Xinit}$ \\
    \ \ \quad where $\|E_{\Yinit,\Xinit}\|_F \le \varepsilon$
    \label{SI_line5}
\STATE Define $\Yhat$ as the leading $\min(|\alpha|,|\beta|)$ 
   columns selected by $P_{\Yinit,\Yinit}$.
\STATE Define $\Xhat$ as the leading $\min(|\alpha|,|\beta|)$ 
   columns selected by $P_{\Xinit,\Xinit}$.
\RETURN Interpolation points $\whX$ and $\whY$
\end{algorithmic}
\end{algorithm}

\ifthenelse{1=1}{}{
\begin{algorithm}
\caption{Skeletonized Interpolation: $\widehat K_{X,Y} = {\sf SI}(X, Y, \varepsilon, {\sf GenInitSet}, r_0)$}
\label{alg:genericSI}
\begin{algorithmic}
\REQUIRE{Kernel $\Kfun$, clusters $X$ and $Y$, tolerance $\varepsilon$, algorithm ${\sf GenInitSet}$, rank $r_0$}
\ENSURE{Approximation $\widehat K_{X,Y}$ of $\KXY$}
\STATE Calculate $(\Xinit, \Wx) = {\sf GenInitSet}(X, r_0)$ and $(\Yinit, \Wy) = {\sf GenInitSet}(Y, r_0)$\\ \label{SI_line1}
\STATE Build
\[ K_w = \Wx^{1/2}\Kinit\Wy^{1/2} \] \label{SI_line2}
\STATE Perform RRQR on $K_w$
\[ K_w P_y = Q_yR_y \] \\    \label{SI_line3}
\STATE Compute the rank $k = \max \{ 1 \le i \le n \; | \;  |R_{y,ii}| \geq \varepsilon |R_{y,11}| \}$ \\
\STATE Define $\Yhat$ as the first $k$ columns selected by the  RRQR pivoting.
\STATE Perform RRQR on $K_w^\top$
\[ K_w^\top P_x = Q_xR_x \]
\STATE Select $\Xhat$ in the same way as $\Yhat$. \label{SI_line6}
\STATE If $|\Xhat| \neq |\Yhat|$, extend the smallest set.
\RETURN $\widehat K_{X,Y} \approx \KXYhat\KXhatYhat^{-1}\KXhatY$
\end{algorithmic}
\end{algorithm}
} 

In this paper, we explore four different heuristics to define those
sets, in which $(\Xinit,\Yinit)$ may be a subset of $(X,Y)$
(endo-skeleton) or not (exo-skeleton). Those four ways are:
\begin{enumerate}
    \item Chebyshev grid (exo)
    \item Random subset of the vertices (endo)
    \item Maximally-Dispersed Vertices  or MDV (endo)
    \item Points on an enclosing surface, e.g., a sphere or ellipsoid (exo).
\end{enumerate}

\subsection{Optimality of the point set}

The question of optimality of the rows and columns used for the pivot block $\KXhatYhat$  has been studied in the past in several papers in a range of fields. See for example~\cite{doi:10.1137/0917055,MIKHALEV2018187,fonarev2016efficient,Schork:2018aa,SAVOSTYANOV2014217,doi:10.1137/1.9781611972825.59,GOREINOV19971,CIVRIL20094801,goreinov1997pseudo,Goreinov2011,tyrtyshnikov1995pseudo}. We summarize some of the key results and their relevance for the current work. In~\cite{doi:10.1137/0917055}, the authors prove that a strong rank-revealing QR (RRQR) factorization can identify $r$ columns such that
\begin{equation*}
  \sigma_i(A \Pi) \ge \frac{\sigma_i(A)}{q(r,n)}
\end{equation*}
for matrix $A$, where $\Pi$ is a matrix with $r$ columns that selects the optimal $r$ columns in $A$, $\sigma_i$ is the $i$th singular value, and $q$ is a ``low-degree'' polynomial. Similarly the error $E$ from a strong RRQR factorization when we keep only the first $r$ columns (rank-$r$ approximation) has an upper bound involving $\sigma_{r+i}$:
\begin{equation*}
  \sigma_i(E) \le \sigma_{r+i}(A) q(r,n)
\end{equation*}
From $A \Pi$, we can similarly calculate a strong RRLQ (the transpose of RRQR) to select $r$ important rows of $A \Pi$. This leads to a square sub-block $A_\si$ of $A$. Note that in this two-step process, the first step with RRQR is the one that determines the accuracy of the low-rank factorization. The second step with RRLQ merely guarantees the stability of the SI procedure (with respect to roundoff errors and small perturbations in the input); specifically, it minimizes the conditioning of the operation with respect to perturbations in the input data. In turn, this guarantees the numerical stability of computing $\KXYhat \KXhatYhat^{-1}$ and $\KXhatYhat^{-1} \KXhatY$, in spite of the fact that the singular values of $\KXhatYhat$ are rapidly decaying.

From the interlacing singular values theorem~\cite{golub13}, the singular values $\{\sigma_1^\si$, \ldots, $\sigma_r^\si\}$, of $A_\si$ keep increasing as $r$ increases. As a result, we have both lower and upper bounds on $\sigma_i^\si$. This leads to a definition of the best choice of SI block. It should satisfy two equivalent properties:
\begin{enumerate}
  \item Its singular values $\{\sigma_i^\si\}_{i = 1,\ldots,r}$ are close to the corresponding singular values of $A$ (i.e., with minimal error).
  \item The volume of $A_\si$, as defined by the absolute value of its determinant, is maximum.
\end{enumerate}
This motivates some heuristic criteria for choosing exo-points:
\begin{enumerate}
  \item They should be easy to calculate, such that computing the exo-points is computationally advantageous compared to selecting endo-points from $\KXY$.
  \item When increasing the size of the matrix, by adding the exo-points $\Xinit$ and $\Yinit$, the increase in rank should be minimal. This happens when the set $(\Xinit,\Yinit)$ is ``close'' to $(X,Y)$; this is important when, for example, the $(X,Y)$ points lie on a sub-manifold.
  \item The volume of $\KXhatYhat$ should be maximum. This implies that the points $(\Xinit,\Yinit)$ are widely spread. In terms of the singular functions of the kernel $\Kfun$, this corresponds heuristically to adding points in regions when the singular functions are rapidly varying (e.g., near the boundary of the domain).
\end{enumerate}
In practice, computing these optimal sets is very expensive since it requires solving an NP-hard combinatorics problem (e.g., finding a permutation that maximizes the volume). But we will see that the heuristics we propose can lead to near-optimal solutions at a fraction of the computational cost.

\subsection{Contribution}

In this paper, we adopted four strategies of computing $\Xinit$ and $\Yinit$ as input for
Algorithm~\ref{alg:genericSI}. Numerical experiments are conducted to compare their performance (in terms of sizes of the initial sets $\Xinit, \Yinit$ and the resulting sets $\Xhat, \Yhat$) of each strategy. The goal is to obtain an algorithm such that
\begin{enumerate}
\item the method has a complexity of $O(r(m+n))$ where $r$ is the target rank;
\item the method is efficient, accurate, and stable;
\item the method is applicable to complicated geometries;
\item the method is simple to implement.
\end{enumerate}

We motivate and explain the four methods of selecting $\Xinit$ 
and $\Yinit$ in \autoref{Initial Interpolation Points Selection}
and
\autoref{Numerical experiments} 
presents numerical experiments to compare the four.

\section{Selecting the initial interpolation points} 
\label{Initial Interpolation Points Selection} 

Consider \autoref{eqn:lowrankKXY} and notice that is can be rewritten
\begin{align*} \KXY \approx \KXYhat \KXhatYhat^{-1} \KXhatY & = (\KXYhat \KXhatYhat^{-1}) \KXhatYhat (\KXhatYhat^{-1} \KXhatY) \\
                                                            & = S_{X,\Xhat} \KXhatYhat T_{Y,\Yhat}^\top \end{align*}
where $S_{X,\Xhat}$ and $T_{Y,\Yhat}$ are Lagrange basis functions (\emph{not} polynomials). Each column can be seen as one of $|\Xhat|$ (resp., $|\Yhat|$) basis functions evaluated at $X$ (resp., $Y$). The multiplication by the node matrix $\KXhatYhat$ produces an interpolator evaluated at $X \times Y$. Because of this representation, the points $\Xhat$ and $\Yhat$ can be considered as interpolation points. Formally, the same statement holds for $\Xinit$ and $\Yinit$, from which $\Xhat$ and $\Yhat$ are constructed.

The selection of the initial sets $\Xinit$ and $\Yinit$ is crucial.  
There are three choices.
\begin{itemize}
\item {\bf Endo-vertices : }
Choose initial points from the inside the set of discretization points,
$\Xinit \subseteq X$.
\item {\bf Exo-vertices : }
Alternatively, choose initial points from outside the set of
discretization points, $\Xinit \cap X = \emptyset$. 
\item {\bf Mixed-vertices : }
If neither $\Xinit \subseteq X$ nor $\Xinit \cap X = \emptyset$ hold,
the $\Xinit$ are mixed, some in $X$, others not in $X$.
\end{itemize}
We have two examples of exo and endo vertices.
\begin{itemize}
\item 
The first exo-method is to choose $\Xinit$ to be a tensor product 
of Chebyshev nodes to enclose the domain $\Xdom$. 
The tensor grid may be skew, and is
chosen to closely enclose the domain.
\item 
The second exo-method is to enclose the domain $\Xdom$ within a
ball centered at the centroid of $\Xdom$, 
and choose $\Xinit$ to be well-dispersed points 
on the surface of the ball.
\end{itemize}
Exo methods have an advantage that once a tensor grid of points 
or points on a sphere have been selected, 
their interpolation matrices can be used for any
domain enclosed by the tensor grid.
This property is not shared by the endo-methods, where the choice of
interpolation points are a function of the domain.
\begin{itemize}
\item 
The first endo-method is to choose 
$\Xinit \subset X$ to be
a set of ``maximally dispersed'' vertices, i.e., we try and maximize
the minimum distance between any two vertices in $\Xinit$.
\item 
The second endo-method is to choose $\Xinit \subset X$ to be 
a random set of vertices.
\end{itemize}
Endo methods have the advantage of being nested, when we need more
accuracy, we can add the next few vertices in the sequence and save
on computation.
\par
We now describe these four methods in detail.
%
%
\ifthenelse{1=1}{}{ 
   \par
   The selection of $\Xinit$ and $\Yinit$ in Skeletonized Interpolation
   is crucial. As described in the introduction, there can be different
   ways to choose $\Xinit$ and $\Yinit$. In this section, we present in
   detail four different ways to build those sets, that is Chebyshev
   nodes, maximally-dispersed vertices (MDV), points on a sphere and
   random vertices.
   
   In particular, one can choose points \emph{within the mesh vertices}
   like MDV or random vertices (``endo" vertices). Or one can choose
   points \emph{outside the mesh vertices}, like Chebyshev nodes or
   points on a sphere (``exo" vertices). In the case of Chebyshev
   nodes, we may need to optimize the construction of the 3D box used
   to construct the nodes (i.e., when the vertices $X$ and $Y$ do not
   lie inside clearly defined boxes), as those sets of points are
   in principle defined over $[-1,1]^d$ (see \autoref{sec:complex}).
   The advantage, however, is that the result is independent from 
   the detail of the location of the mesh vertices.
} 
%
%
\par
\subsection{Tensor grids formed of Chebyshev nodes} \label{sec:cheb}
\par
Polynomial interpolation at Chebyshev nodes is a popular approach to
approximate low-rank kernel matrices, used for time-domain cosine 
expansion of large stationary covariance matrices \cite{van2014}, 
and general least-square analysis \cite{Coles2011}.
\par
The idea of Chebyshev interpolation~\cite{Fong2009TheBF,Messner2012},
is simple ---
the more differentiable the kernel
function $\Kfun$, the faster the coefficients of its Chebyshev
expansion decay.  
This means that the kernel matrix has singular values that decay,
and a low rank representation of the matrix is a good approximation.
This is independent of the choice of discretization points $X$ and
$Y$ inside domains $\Xdom$ and $\Ydom$.
\par
Consider a kernel $\Kfun$ defined over the domains 
$\Xdom \subseteq {\mathbb R}^d$ and $\Ydom \subseteq {\mathbb R}^d$. 
The kernel can be approximated using a polynomial interpolation rule as 
\begin{align} \label{eq:chebyshev_interpolation}
\Kfun(x,y) 
&\approx \sum_{x^{\circ} \in X^{\circ}} \sum_{y^{\circ} \in Y^{\circ}}
S_{x^{\circ}}(x)\,\Kfun(x^{\circ},y^{\circ})\, T_{y^{\circ}}(y)
\end{align}
where $\Xinit$ and $\Yinit$ 
are Chebyshev interpolation points
and $S_{x^{\circ}}(x)$ and $T_{y^{\circ}}(y)$ are Lagrange polynomials. In 1D over $[-1, 1]$, the $m$ Chebyshev
nodes (and their associated integration weights) are defined as
\[ x^\circ_k = \cos\left( \frac{2k-1}{2m}\pi \right), \quad w^\circ_k = \frac{\pi}{m} \sin\left( \frac{2k-1}{2m} \pi \right) \text{ for } k = 1,\dots,m. \]
See \cite{BerrutTrefethen2004} for further details about practical Chebyshev interpolation.
\par
In order to accomodate arbitrary geometries, we build a tensor grid where the numbers of points along each 
dimension are $n_1$, $n_2$ and $n_3$ for a total number of points
$N = n_1 n_2 n_3$. 
The $n_i$ are chosen in rough proportion to the lengths of the sides
of the enclosing box.
We use a simple principle component analysis (PCA) to
find the orientation of the bounding box.
The lengths of the sides of the box are chosen to enclose the domain.
This defines $\Xinit$ and $\Yinit$.
The associated weights $\Wx, \Wy$ are the products of the associated one-dimensional integration weights.
\par
Note that while this approach is justified by the existance of the interpolation (cf.~\autoref{eq:chebyshev_interpolation}), 
we merely rely on the nodes and the weights. There is no need to ever build or evaluate the associated Lagrange basis functions.
\par
Benchmark tests in \cite{LeopoldSI2017} show that
{\sf SI-Chebyshev} works well when the sizes of $X$ and $Y$ are
large, $O(10,000)$. 
However, when $X$ and $Y$ are small, $O(100)$ , the construction 
of the interpolation grids from Chebyshev expansion may be inefficient 
(the initial rank $r_0$ has to be too large for a given tolerance).
\par
Since the 1-D Chebyshev points are clustered toward the endpoints
of the interval, the tensor grid points are clustered in the corners.
Unless the domains fill out the corners of the enclosing boxes,
a corner may not be the best place to place interpolation points.
The second exo-method distributes the points in a more even fashion.
\par
\subsection{Points on a sphere} \label{sec:sphere}
\par
The second method of constructing exo-vertices $\Xinit$ 
and $\Yinit$ is to evenly distribute points on a sphere 
that encloses the domain. 
For kernels that satisfy Green's theorem, points on an
enclosing surface are sufficient. 
Since for a geometry of dimension $d$, we only need to build 
$\Xinit$ and $\Yinit$ on a manifold of dimension $d-1$, this 
approach can significantly reduce the size of $\Xinit$ and $\Yinit$.
\par
Green's third identity represents the far-field using an 
integration over the boundary instead of the whole domain. 
Consider points $x$ in $\Xdom$ and $y$ in $\Ydom$. 
Define a surface $\Gamma$ enclosing $x$. 
The exterior domain
$\Omega \subset {\mathbb R}^3$ with boundary $\Gamma$ contains $y$
but not $x$.
\par
The function $\Kfun(x, y)$ satisfies the boundary value problem
\begin{alignat*}{2}
  \Delta_y u(y) & = 0 && y \in \Omega \\
  u(y) & = \Kfun(x, y) \quad && \text{for all $y \in \Gamma$}
\end{alignat*}
Generally, the function $u$ satisfies the following representation
formula~\cite{ banerjee1981boundary, becker1992boundary,
               brebbia1980boundary, hall1994boundary}
\begin{equation} \label{eq:greenDeri}
  u(y) = \int_\Gamma \left\{ [u(y)]_\Gamma \; \frac{\partial \Kfun}{\partial n_1}(z,y)
- \Big[\frac{\partial u}{\partial n_2}(x,z) \Big]_\Gamma
\Kfun(z,y) \right\} \text{d}z
\end{equation}
for all $y \in {\mathbb R}^3 \setminus \Gamma$, where $[\;]_\Gamma$
denotes the jump across $\Gamma$. 
Assume the field $u$ is continuous across the boundary $\Gamma$. 
Note that
it is then equal to $\Kfun(x, y)$ for $y \in \Omega$, but, being
smooth for all $y \in {\mathbb R}^3 \setminus \Gamma$, $u$ is not
equal to $\Kfun(x, y)$ for $y \not\in \Omega$ (for example near $x$
where $\Kfun(x, y)$ is singular). 
With this choice $[u(y)]_\Gamma = 0$, 
we simplify \autoref{eq:greenDeri}.
\begin{equation} \label{eq:green}
\Kfun(x, y) = - \int_\Gamma 
  \Big[\frac{\partial u}{\partial n_2}(x,z) \Big]_\Gamma \;
\Kfun(z,y) \; \text{d}z
\end{equation}
This implies that it is possible to find ``pseudo-sources'' on the
surface of $\Gamma$ that will reproduce the field $\Kfun(x, y)$ on
$\Omega$~\cite{Makino1999}. This is the motivation behind using
points on a sphere, since for instance $\Kfun(x, y) =
\|x-y\|_2^{-1}$ satisfies the potential equation. However, not all
kernels satisfy this equation. In particular $\Delta_x
\|x-y\|_2^{-2} \neq 0$. As a result, one can't apply Green's
theorem, and points on a sphere are not enough to interpolate this
kernel. 
\par
This is illustrated in \autoref{plate-coil computational results}.
A representation similar to \autoref{eq:greenDeri} has been used before.
One idea, the Green hybrid method (GrH), is explored in
\cite{Borm2014GreenHybrid}. This method takes advantage of a
two-step strategy. It first analytically approximates the kernel
using an integral representation formula, then further compresses
its rank by a factor of two using a nested cross approximation. 
Another approach is to use Green's formula to place a bounding
circle around the domain, and then spread interpolation points 
equally spaced around the circle.
\par
In \cite{Borm20132DSphere} and \cite{Borm2014GreenHybrid},
both $\Kfun$ and its normal derivative appear, 
while our approach only uses $\Kfun$ on the uniformly distributed
points on the bounding sphere.
Similar to our method, ``pseudo-points'' \cite{Makino1999}
anchor the multipole expansion to approximate the potential field. 
\par
In this paper, to deal with 3D geometries, 
we use a Fibonacci lattice
\cite{gonzalez2010measurement,hardy1979introduction} 
to spread uniformly distributed points ${\bf x} \in \Xinit$ 
on a 3D bounding sphere.
\begin{center}
\includegraphics[width=.25\textwidth]{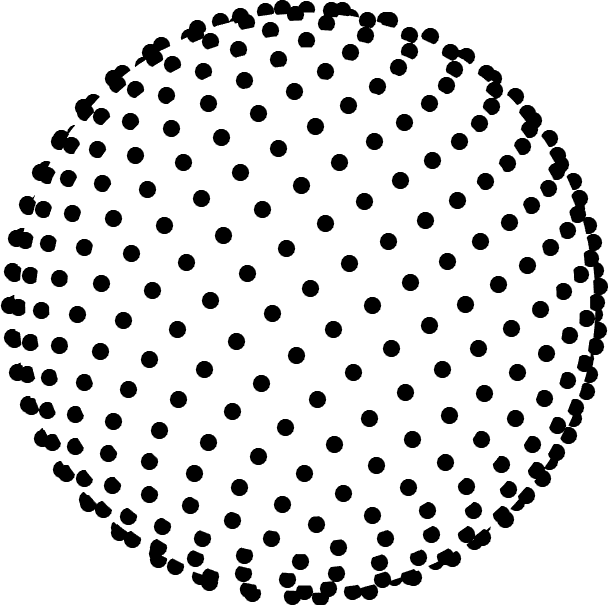}
\raisebox{1.5cm}{
\begin{minipage}{2.8 in}
\begin{equation*}
{\bf x}_i = c_M + r_M
\begin{bmatrix} r_i \cos\theta_i &
r_i \sin\theta_i &
z_i \end{bmatrix}^\top
\end{equation*}
\begin{align*}
\Delta \theta &= (3-\sqrt{5}) \pi & \Delta z &= \frac{2}{n-1} \\
\theta_1 &= \Delta \theta & z_1 &= -1 \\
\theta_{i+1} &= \theta_i + \Delta \theta 
& z_{i+1} &= z_i + \Delta z \\
r_i & = \sqrt{1 - z_i^2} 
\end{align*}
\end{minipage}
}
\end{center}
An obvious extension is to include points interior to the sphere,
i.e., place interpolation points on nested spheres.
\par
\subsection{Random vertices}
\label{subsec:random}
\par
The first endo-method is very simple, simply select a set of
interpolation vertices $\Xinit$ from $X$ at random
\cite{halko2011finding,
      gu2016efficient,
      martinsson2016randomized,
      mary2015performance}.
For large $X$, this is reasonable, but in general,
this approach is not as good as our second endo-method.
\par
The scaling matrix $W_{\Xinit,\Xinit}$ is diagonal,
with $W_{x_i,x_i}$ proportional to the ``area'' of a patch
surrounding vertex $x_i$. We can define area this patch 
as the weight of the vertices that are close to $x_i$.
\par
Consider a vertex $v \in X$. 
If there is one interpolation vertex $x_i \in \Xinit$ that 
it is closest to, then give $x_i$ the weight of vertex $v$.
If there are two or more interpolation points that are closest to
$v$, then give an equal portion of the weight of $v$ to each of
the closest interpolation points.
With this definition of area, we define the scaling matrices
$W_{\Xinit,\Xinit}$ and $W_{\Yinit,\Yinit}$.
\par
\subsection{Maximally dispersed vertices}
\label{subsec:MDV}
\par
If domains $\Xdom$ and $\Ydom$ are well-separated,
the range of $\KXY$ is well represented by interpolation vertices 
$\Xinit$ and $\Yinit$ that are roughly uniformly distributed 
throughout $X$ and $Y$.  
\par
We construct a sequence of vertices $x_1, x_2, \ldots, x_m$ 
such that each leading subset is ``maximally dispersed'',
they are as far away from each other as possible.
\par
We want to choose interpolation vertices $\Xinit$
from discretization vertices $X$.
Choose a random vertex $u$ and find a vertex $v$ the furthest from $u$.
The first interpolation vertex in $\Xinit$ is vertex $v$.
To add another interpolation vertex, we look for a vertex that is
the furthest distance from any interpolation vertex, 
choose $v \in X \setminus \Xinit$ 
with a maximum minimum distance.
\begin{equation}
\label{eqn:MDV}
\displaystyle
\min_{x^{\circ} \in \Xinit} {\sf dist}(v,x^{\circ})
=
\max_{w \in X \setminus \Xinit} 
\left(
\min_{x^{\circ} \in \Xinit} {\sf dist}(w,x^{\circ})
\right)
\end{equation}
We use either the Euclidean distance metric 
or the graph distance metric.
When we build the interpolation set up one vertex at a time,
and at each step the equality holds, we say that these vertices
are {\bf maximally dispersed}.
\begin{figure}[!ht]%
\def\la{0.23\textwidth}%
\centering%
\subfloat[Select 1 MDV \label{fig:MDV1}]{
    \includegraphics[width=\la]{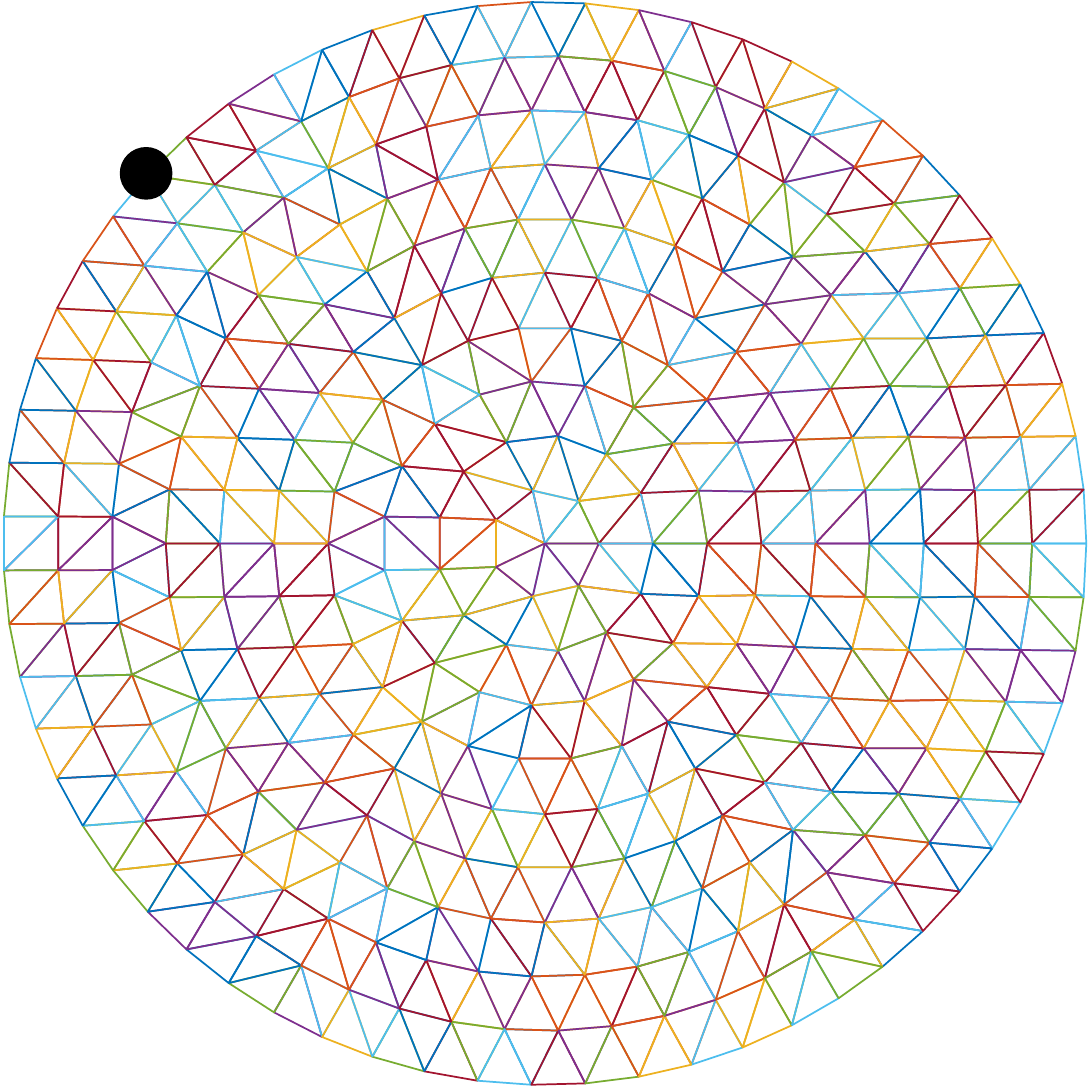}
}\hfill%
\subfloat[Select 10 MDV \label{fig:MDV10}]{
    \includegraphics[width=\la]{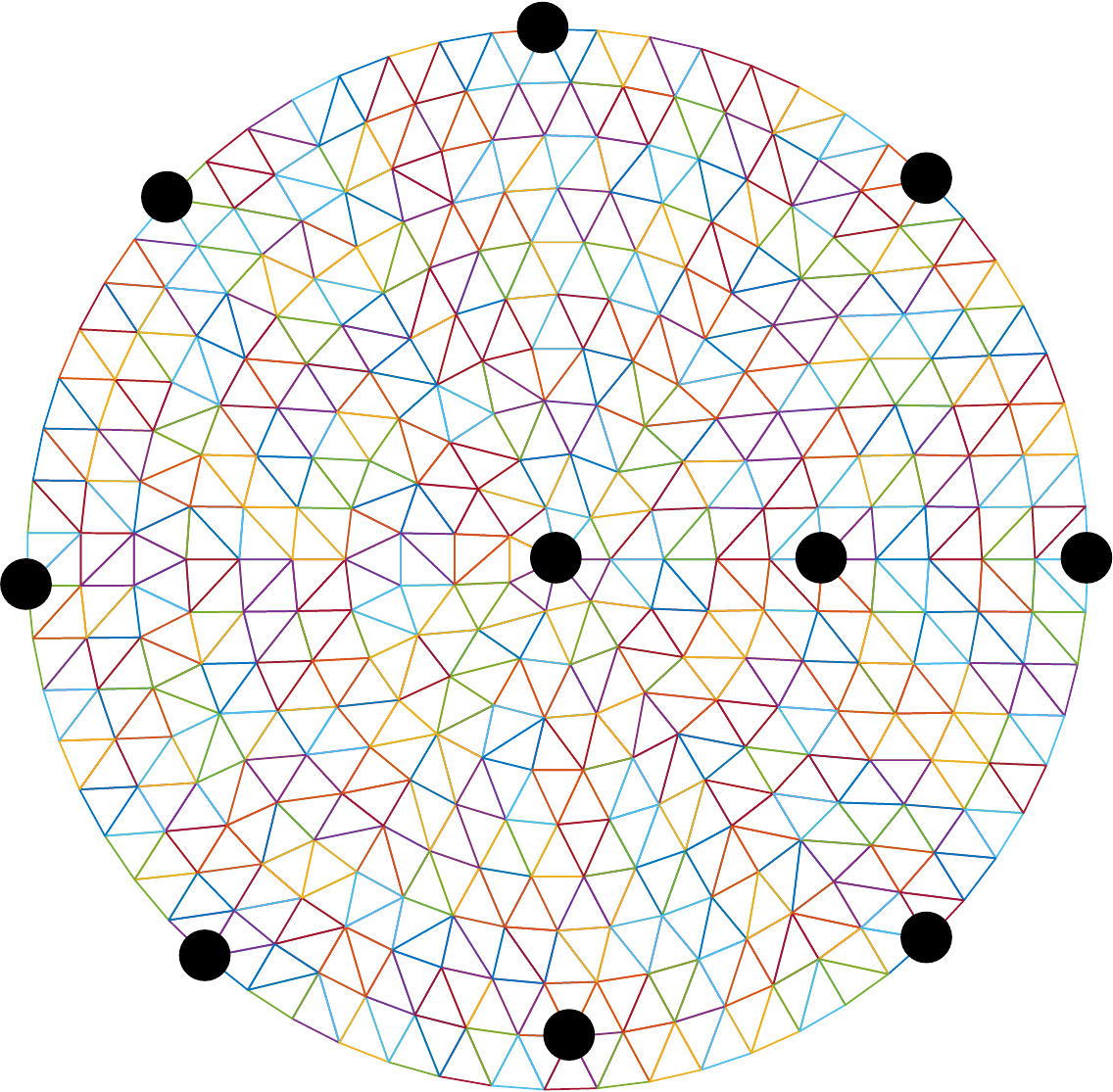}
}\hfill%
\subfloat[Select 25 MDV \label{fig:MDV25}]{
    \includegraphics[width=\la]{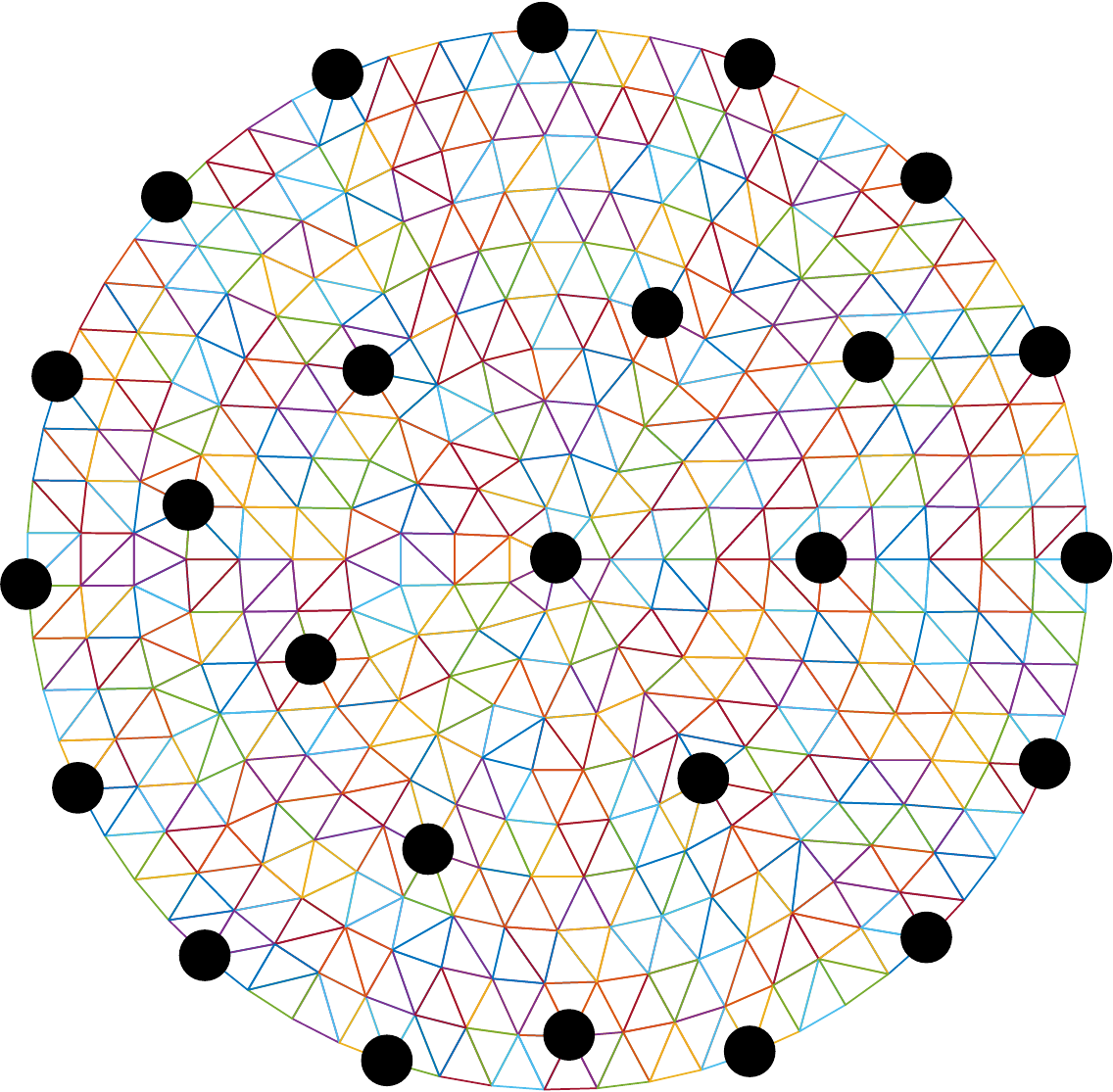}
}\hfill%
\subfloat[Select 50 MDV \label{fig:MDV50}]{
    \includegraphics[width=\la]{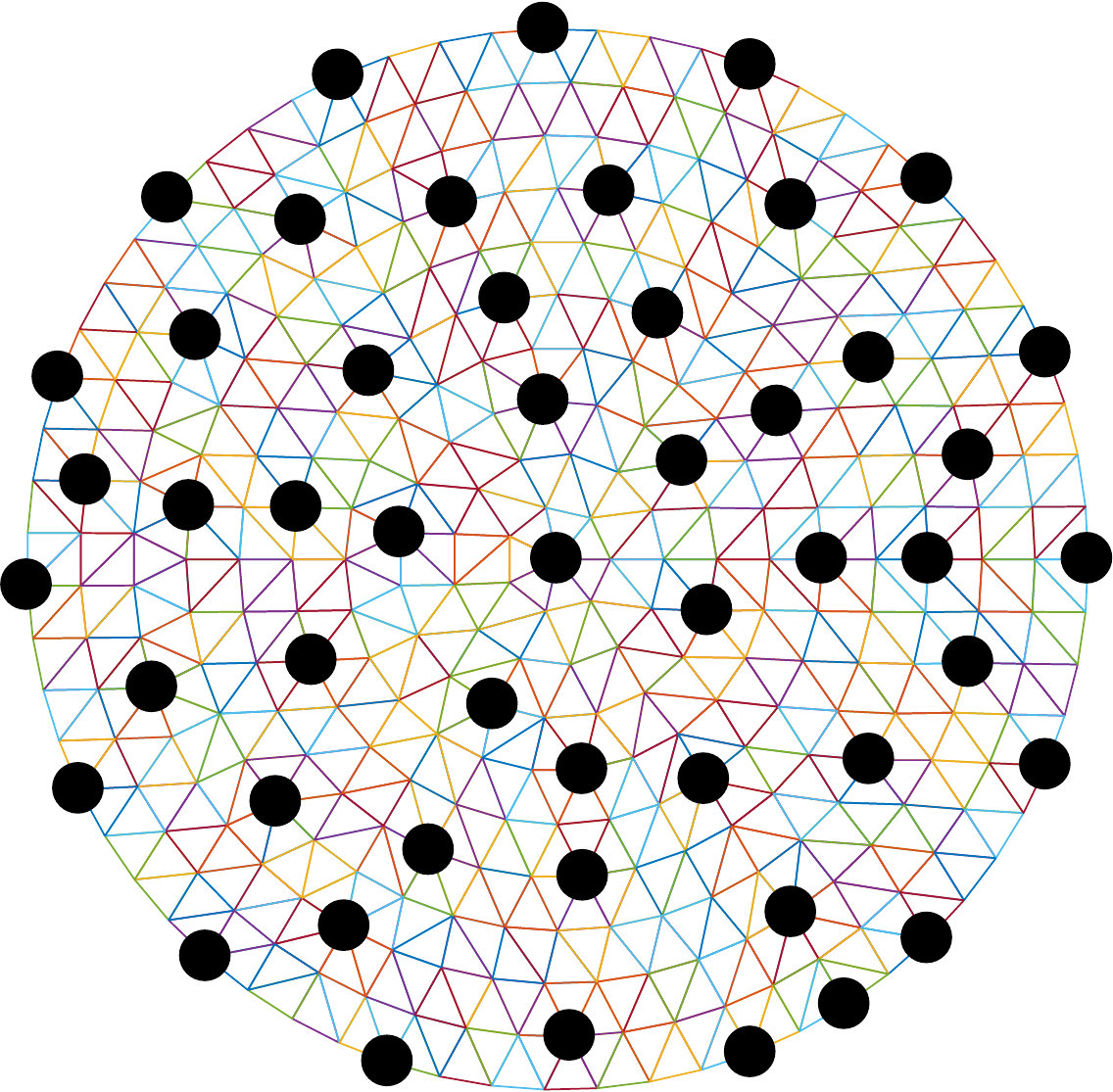}
}%
\caption{MDV initial interpolation points $\Xinit$ 
         for a circular disk $\Xdom$,
         triangulated with 351 discretization points $X$.}%
\label{fig:MDVinAction}%
\end{figure}
\autoref{fig:MDVinAction} shows the process for 1, 10, 25 and 50
vertices in $\Xinit$ for a triangulated disk with $|X| = 351$.
The scaling matrices for MDV are area-based, described in
\autoref{subsec:random}.

\ifthenelse{1=1}{}{ 

With this, one can then use \autoref{alg:genericSI} with 
${\sf GenInitSet}(X,r) = (\Xinit, \Wx)$ defined by 
$\Xinit = {\sf SI-MDV}(X,r,\dist)$ and $\Wx$ the identity. 
This defines the {\sf SI-MDV} algorithm.

} 

\subsection{Discussion}
\label{subsec:Discussion}
We compare four methods :
exo-methods {\sf SI-Chebyshev} and {\sf SI-sphere}, and
endo-methods {\sf SI-MDV} and {\sf SI-random}.
    \begin{figure}

        \centering
        \begin{tikzpicture}
            \begin{semilogyaxis}[
            width=8cm, 
            height=6cm,
            xlabel={Rank},
            ylabel={Relative Frobenius Error},
            grid = major,
            legend entries={No weights, Weights},
            ]
            \addplot [red,mark=diamond*,mark size=1pt] table [x=rank, y=no_weights] {Images/weights_no_weights.dat};
            \addplot [blue,mark=square*,mark size=0.5pt] table [x=rank, y=with_weights] {Images/weights_no_weights.dat};
            \end{semilogyaxis}
        \end{tikzpicture}
        \caption{(Relative) Frobenius error $\|\KXY - \KXYhat \KXhatYhat^{-1} \KXhatY\|_F$ when $\Xhat$ and $\Yhat$ are built using \autoref{alg:genericSI} (with the rank fixed a priori) with and without the weight matrices. $X$ and $Y$ are $20^3$ points on two facing unit-cubes separated by a distance of 1, and $\Xinit$, $\Yinit$ are tensor grids of $7^3$ Chebyshev nodes. $\Kfun(x, y) = 1/r$. }
        \label{fig:weights_no_weights}
    \end{figure}
\begin{itemize}
\item {\bf Scaling matrices} 
    The weight matrices $\Wx$ and $\Wy$ are needed so that the 2-norm of the rows (resp. columns) of $\widehat K_{\Xinit, \Yinit}$ properly approximates the $L_2$ norm of $\Kfun$ over $\Xdom$ (resp. $\Ydom$). Otherwise, a higher concentration of points in a given area may be given a higher than necessary importance during the columns (resp. rows) selection process in \autoref{alg:genericSI}. Figure \ref{fig:weights_no_weights} illustrates the impact on the Frobenius error (a proxy for the $L_2$ error over $\Xdom \times \Ydom$) when using weights or not with {\sf SI-Chebyshev}. We observe that the absence of weights leads to a higher Frobenius error (with a higher variance) in the result.
    This justifies the choice of weights matrices for the four methods:
\begin{itemize}
\item[$\bullet$]
{\sf SI-sphere} does a very good job of distributing the points
evenly on the sphere, so the areas of the patches are nearly constant. 
The scaling matrix for the unit sphere 
is close to $(4\pi/|\Xinit|)$ times the identity.
\item[$\bullet$]
{\sf SI-MDV} also does a good job of evenly distributing vertices,
see \autoref{fig:MDV50}. 
Its scaling matrix is nearly a constant diagonal.
\item[$\bullet$]
For moderate numbers of initial vertices, {\sf SI-random} can have
some variation in areas, and so a non-unit scaling matrix should be
used.
\item[$\bullet$]
The nodes in a Chebyshev tensor grid are found mostly near the
boundaries, the edges, and the corners. The scaling matrix 
for {\sf SI-Chebyshev} is necessary for the selection of 
interpolation points $\whX$.
\end{itemize}
\item {\bf Directional bias} 
\par
All methods that evenly distribute points can suffer from loss of
accuracy for near-field matrices. 
Consider below where we show four steps of choosing a good set of
dispersed vertices for a rectangular domain $X$.
The $Y$ domain is two diameters to the northeast of $X$. 
We compute the RRQR factorization of the kernel matrix $\KXY^T$ 
and show the first 5, 10, 15 and 20 vertices chosen as pivots.
There is a definite bias to the northeast, in the direction of
the $Y$ domain. When the two domains are closer, this bias is more
pronounced.
\begin{minipage}{4.5 in}
\centering
\includegraphics[width=0.2\textwidth]{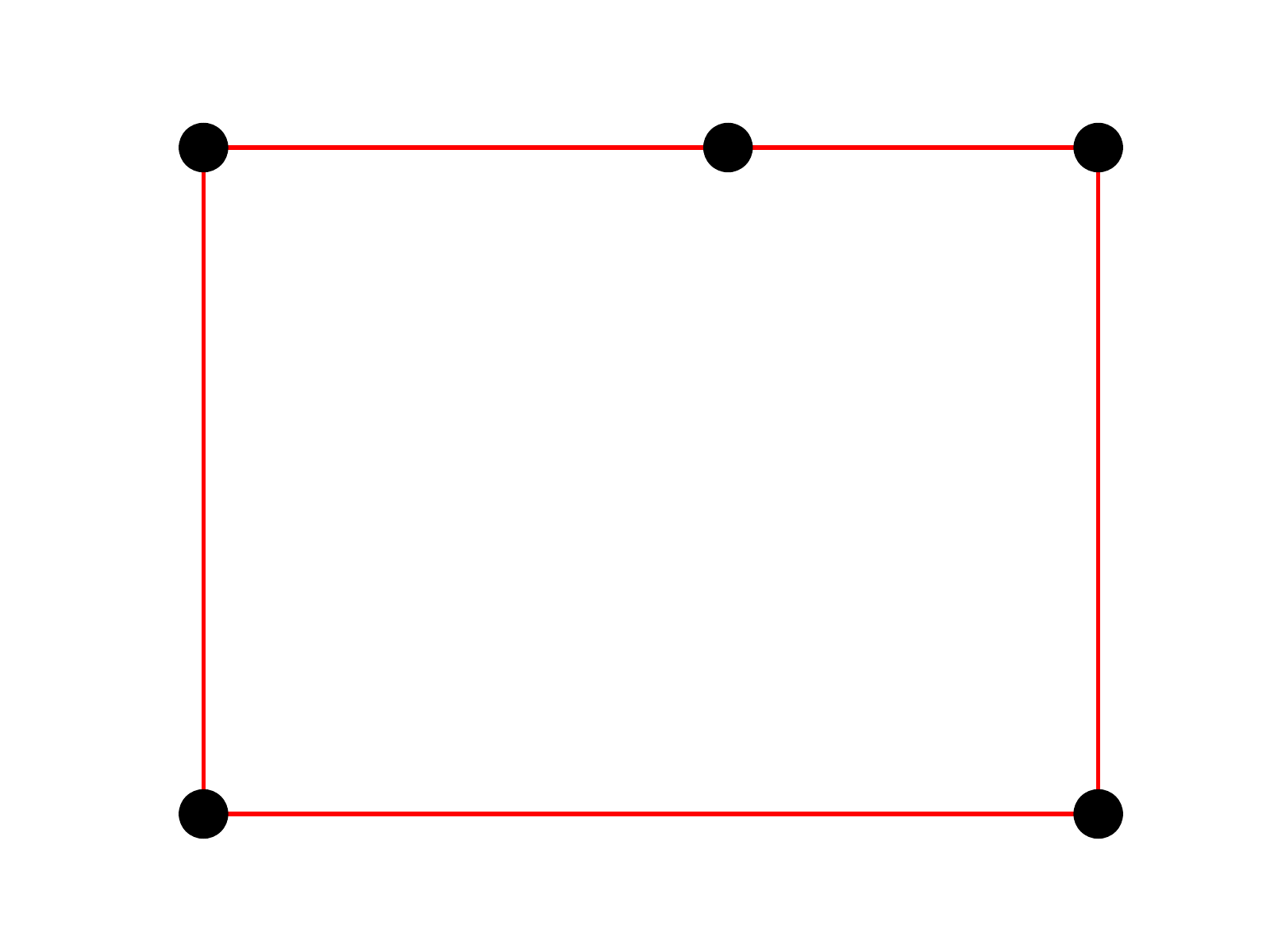}
\includegraphics[width=0.2\textwidth]{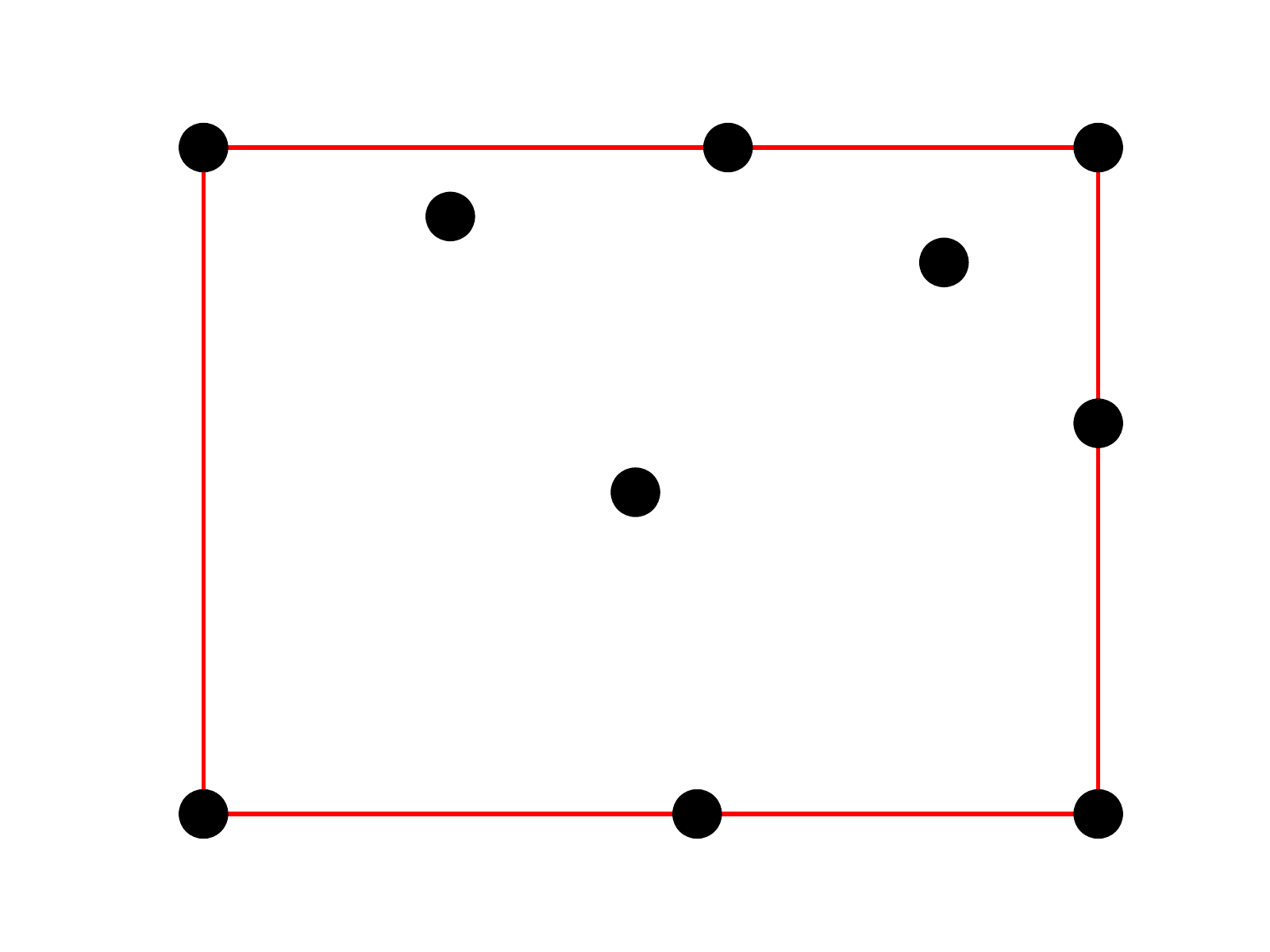}
\includegraphics[width=0.2\textwidth]{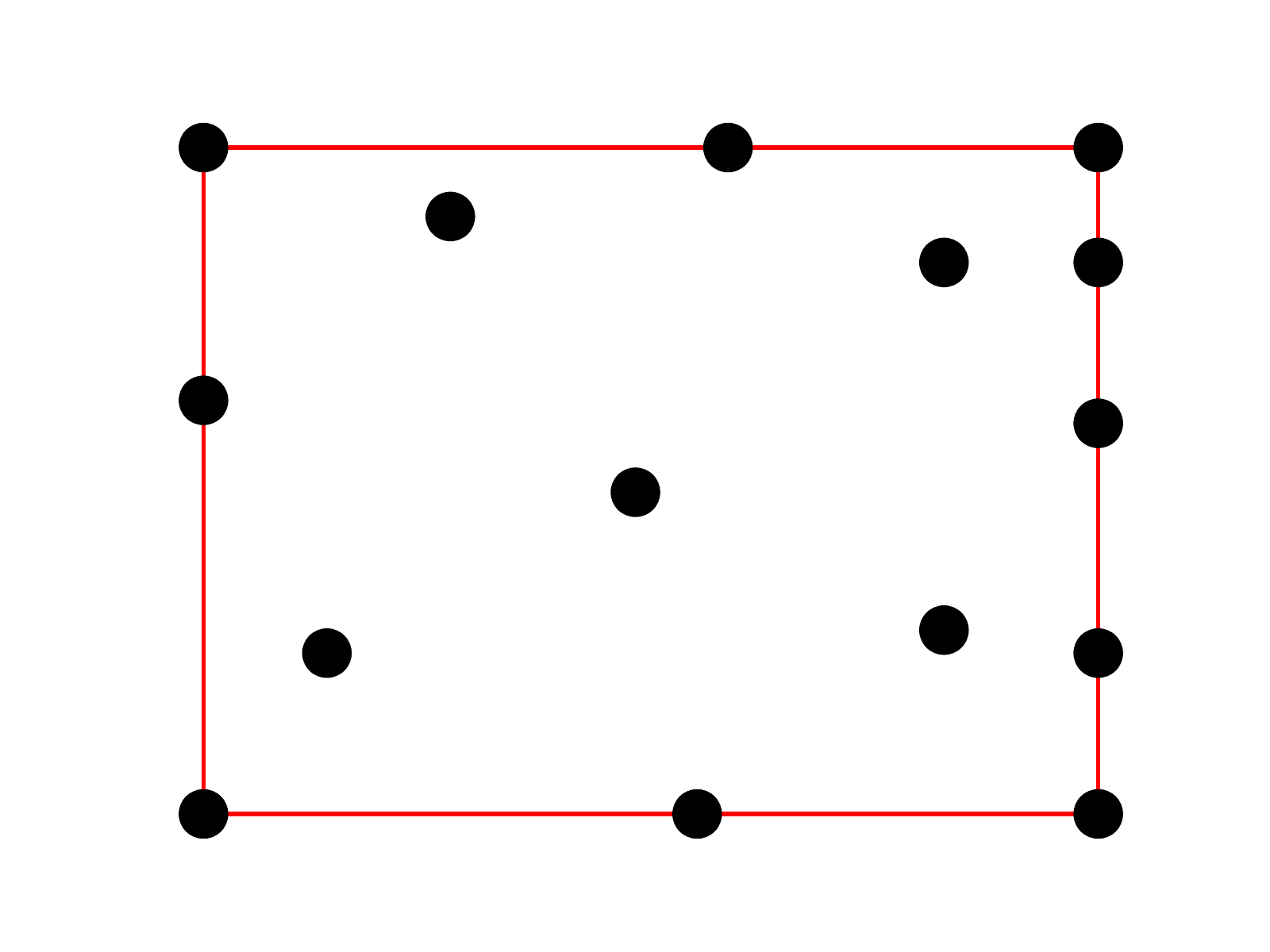}
\includegraphics[width=0.2\textwidth]{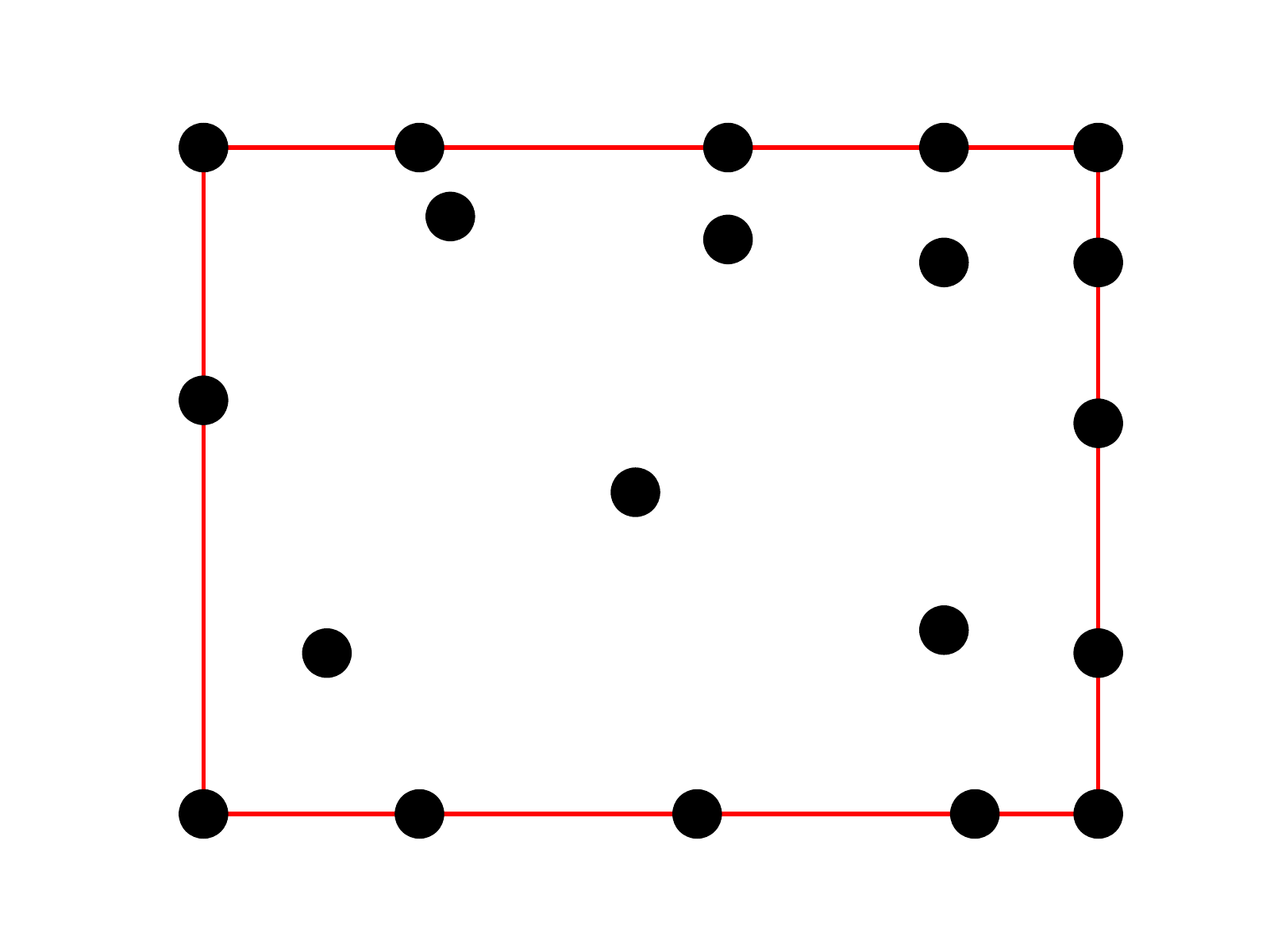}
\end{minipage}
\par
{\sf SI-MDV}, {\sf SI-sphere} and so some extent {\sf SI-random} will 
need a more dense uniformly dispersed initial set $\Xinit$ to be able 
to provide enough vertices in the necessary area of concentration.
Conversely, {\sf SI-Chebyshev} may not suffer too badly, since its
nodes are largely found near the boundary.
\item{\bf Further compression of the low-rank approximation} 
\par
Given an approximation 
\begin{equation}
\| \KXY - \KXYhat \KXhatYhat \KXhatY \|_F \leq \varepsilon \| \KXY \|_F
\end{equation}
one can \emph{always}, no matter the sets $\Xhat$ and $\Yhat$,
obtain a re-compressed approximation with a rank equal to or
near the optimal $\varepsilon$-rank of $\KXY$.
\par
To see this, consider a matrix $A$ and a matrix $B$ of rank $r_1$
(in our problem, $A = \KXY$ and $B = \KXYhat \KXhatYhat^{-1}
 \KXhatY$). The key is that 
\begin{equation}
\| A - B \|_2 \leq \varepsilon 
\Rightarrow | \sigma_k(A) - \sigma_k(B) | \leq \varepsilon. 
\end{equation}
This is a direct consequence of corollary 8.6.2 in \cite{golub13}.
Intuitively, this means the singular values are continuous, and so up 
to $\varepsilon$, $A$ and $B$ have roughly the same singular values. 
One can compress $B$ to another low-rank matrix with an even smaller 
rank, $r_2$, closer to the $\varepsilon$-rank of $A$.
\begin{align} 
\notag
\KXYhat \KXhatYhat^{-1} \KXhatY 
&= \left( Q_{X,\alpha} R_{\alpha,\whY} \right) \KXhatYhat^{-1} 
   \left(Q_{Y,\beta} R_{\beta,\whX} \right)^\top 
\\ \notag
&= Q_{X,\alpha} 
   \left( R_{\alpha,\whY} \KXhatYhat^{-1} R_{\beta,\whX}^\top \right) 
   Q_{Y,\beta}^\top 
\\ \notag
&\approx Q_{X,\alpha} 
   \left( U_{\alpha,\gamma} \Sigma_{\gamma,\gamma} 
          V_{\beta,\gamma}^\top \right) 
   Q_{Y,\beta}^\top 
\\ \label{eqn:compress}
&= \left( Q_{X,\alpha} U_{\alpha,\gamma} \right)
   \left( \Sigma_{\gamma,\gamma} Q_{Y,\beta} V_{\beta,\gamma} \right) 
 = U_{X,\gamma} W_{Y,\gamma} 
\end{align}
The error in the truncation of the singular values of $R_{\alpha,\whY} \KXhatYhat^{-1} R_{\beta,\whX}^\top$ (from the second to third line) has not been further 
amplified since $Q_{X,\alpha}$ and $Q_{Y,\beta}$ are orthogonal.

\ifthenelse{1=1}{}{
\begin{algorithm}
\caption{Recompression : $UV^\top = {\sf Recomp}(\Kfun, X, Y, \Xhat, \Yhat, \varepsilon)$}
\label{alg:recomp}
\begin{algorithmic}
\REQUIRE{Approximation $\KXY \approx \KXYhat \KXhatYhat^{-1} \KXhatY$ of rank $r_1$}
\ENSURE{Recompression $\KXY \approx UV^\top$ of rank $r_2 \leq r_1$}
\STATE Compute the QR factorizations
\[ Q_x R_x = \KXYhat, \quad
Q_y R_y = \KXhatY^\top \]

Compute the SVD
\[ U_x S U_y = R_x \KXhatYhat^{-1} R_y^\top \]

Truncate $S$ at $\varepsilon$ and obtain the rank $r_2$
\[ r_2 = \max\{1 \leq i \leq r_1 | S_{11} \varepsilon \leq S_{ii} \} \]

\RETURN $UV^\top$ with $U = Q_x U_{x|:,1:r_2}$ and $V = Q_y U_{y|:,1:r_2} S_{1:r_2,1:r_2}$
\end{algorithmic}
\end{algorithm}
} 
The singular values of $\KXY$ 
and $R_{\alpha,\whY} \KXhatYhat^{-1} R_{\beta,\whX}^\top$
are close up to $\varepsilon$. 
Hence, by truncating up to $\varepsilon$, 
we can expect to recover the $\varepsilon$ rank of $\KXY$.
We illustrate this result on a concrete example in 
\autoref{sec:comp_res}.
\item{\bf Computational complexities}
\par
We summarize the
computational complexities of the Skeletonized Interpolation and
the recompression process, where we assume 
\begin{equation*}
|X|=|Y|=n,\ 
|\Xinit|=|\Yinit| = r_0 \le n, \ 
\text{\ and\ \ }
|\whX|=|\whY| = r_1 \le r_0.
\end{equation*}
Furthermore, we assume there that the cost of
evaluating the kernel on $X \times Y$ is $O(|X| |Y|)$.
\par
For Skeletonized Interpolation, the leading cost is the 
{\sf GenInitSet} algorithm (if using MDV), the RRQR's and the
construction of the left and right factors. 
Many algorithms exist to compute RRQR's; in our experiments, 
we use a simple column-pivoted QR algorithm with cost $O(r_0^2 r_1)$
(\cite{golub13}, page 278).
\par \medskip \par
\begin{center}
\begin{tabular}{l|c|c}
{\bf Skeletonized Interpolation}
   & kernel & linear \\ 
(\autoref{alg:genericSI}) & evaluation & algebra \\ \hline
{\sf GenInitSet} & & $O(nr_0)$ for MDV \\
                 & & $O(r_0)$ otherwise \\
build $K_{\Xinit,\Yinit}$ & $O(r_0^2)$ & \\
RRQR of $K_{\Xinit,\Yinit}$ $K_{\Xinit,\Yinit}^\top$ 
        & & $O(r_0^2r_1)$ \\
build $K_{\Xinit,\whY}$ and $K_{\Yinit,\whX}$ & $O(nr_1)$ & \\
build $K_{\whX,\whY}$ & $O(r_1^2)$ & \\
LU factorization of $K_{\whX,\whY}$ & & $O(r_1^3)$ \\ \hline
total & $O(nr_1 + r_0^2)$ & $O(nr_0 + r_0^2r_1)$ for MDV \\
      &                   & $O(r_0^2r_1)$ otherwise \\
\end{tabular}
\end{center}
The goal is to have $r_1$ be close to the numerical rank of $\KXY$,
but if it is larger, as we will see is the case for {\sf SI-Chebyshev},
then it will pay to recompress the matrix.
\par \medskip \par
\begin{center}
\begin{tabular}{l|c}
{\bf Recompression} \  (\autoref{eqn:compress}) & linear algebra 
   \\ \hline
QR of $K_{X,\whY}$ and $K_{\whX,Y}^\top$ & $O(nr_1^2)$ \\
compute $R_{\alpha,\whY} K_{\whX,\whY}^{-1} R_{\beta,\whX}^T$
      & $O(r_1^3)$. \\
SVD factorization of 
      $R_{\alpha,\whY} K_{\whX,\whY}^{-1} R_{\beta,\whX}^T$
      & $O(r_1^3)$ \\ 
compute $U_{X,\gamma}$ and $V_{Y,\gamma}$ & $O(nr_1r_2)$ \\ \hline
total & $O(nr_1^2)$ 
\end{tabular}
\end{center}

\ifthenelse{1=1}{}{ 


The second method of constructing exo-vertices as $\Xinit$ and $\Yinit$ is building uniformly
distributed points on a bounding sphere. For kernels that satisfy Green's theorem, points on an enclosing surface are sufficient. Since for a geometry of dimension $d$, we only need to build $\Xinit$ and $\Yinit$ on a manifold of dimension $d-1$, this approach reduces the size of $\Xinit$ and $\Yinit$ significantly.

Mathematically, Green's third identity can be used to represent the far-field using an integration over the boundary instead of the whole domain. This can be derived as follows. Consider a point $x$ in $\mathcal{X}$ and $y$ in $\mathcal{Y}$. Define a surface $\Gamma$ enclosing $x$. We define $\Omega \in {\mathbb R}^3$ with boundary $\Gamma$, containing $y$; $\Omega$ is the exterior domain extending to infinity; we have $x \in {\mathbb R}^3 \setminus \Omega$. The function $\Kfun(x, y)$ satisfies the boundary value problem
\begin{alignat*}{2}
  \Delta_y u(y) & = 0 && y \in \Omega \\
  u(y) & = \Kfun(x, y) \quad && \text{for all $y \in \Gamma$}
\end{alignat*}
Generally, the function $u$ satisfies the following representation formula~\cite{hall1994boundary,banerjee1981boundary,brebbia1980boundary,becker1992boundary}
\begin{equation} 
  u(y) = \int_\Gamma \left\{ [u(y)]_\Gamma \; \frac{\partial \Kfun}{\partial n_1}(z,y)
- \Big[\frac{\partial u}{\partial n_2}(x,z) \Big]_\Gamma
\Kfun(z,y) \right\} \text{d}z
\end{equation}
for all $y \in {\mathbb R}^3 \setminus \Gamma$, where $[\;]_\Gamma$
denotes the jump across $\Gamma$. Let us assume that the field $u$
is chosen to be continuous across the boundary $\Gamma$. Note that
it is then equal to $\Kfun(x, y)$ for $y \in \Omega$, but, being
smooth for all $y \in {\mathbb R}^3 \setminus \Gamma$, $u$ is not
equal to $\Kfun(x, y)$ for $y \not\in \Omega$ (for example near $x$
where $\Kfun(x, y)$ is singular). With this choice
$[u(y)]_\Gamma = 0$, and we get:
\begin{equation} 
  \Kfun(x, y) = - \int_\Gamma \Big[\frac{\partial u}{\partial n_2}(x,z) \Big]_\Gamma \;
\Kfun(z,y) \; \text{d}z
\end{equation}
This implies that it is possible to find ``pseudo-sources'' on the
surface of $\Gamma$ that will reproduce the field $\Kfun(x, y)$ on
$\Omega$~\cite{Makino1999}. This is the motivation behind using
points on a sphere, since for instance $\Kfun(x, y) =
\|x-y\|_2^{-1}$ satisfies the potential equation. However, not all
kernels satisfy this equation. In particular $\Delta_x
\|x-y\|_2^{-2} \neq 0$. As a result, one can't apply Green's
theorem, and points on a sphere are not enough to interpolate this
kernel. 
This is illustrated in \autoref{plate-coil computational results}.
\par
A representation similar to \autoref{eq:greenDeri} has been used before.
One idea, the Green hybrid method (GrH), is explored in
\cite{Borm2014GreenHybrid}. This method takes advantage of a
two-step strategy. It first analytically approximates the kernel
using an integral representation formula, then further compresses
its rank by a factor of two using a nested cross approximation. 
Another approach is to use Green's formula to place a bounding
circle around the domain, and then spread interpolation points 
equally spaced around the circle.
\par
In \cite{Borm20132DSphere} and \cite{Borm2014GreenHybrid},
both $\Kfun$ and its normal derivative appear, 
while our approach only uses $\Kfun$ on the uniformly distributed
points on the bounding sphere.
Similar to our method, ``pseudo-points'' \cite{Makino1999}
anchor the multipole expansion to approximate the potential field. 
\par
In this paper, to deal with 3D geometries, we construct 3D bounding spheres and spread uniformly distributed points as $\Xinit$ and $\Yinit$. The generation of such points is introduced in \autoref{SpherePointsGeneration}.



This implies that it is possible to find ``pseudo-sources'' on the surface of $\Gamma$ that will reproduce the field $\Kfun(x, y)$ on $\Omega$~\cite{Makino1999}. This is the motivation behind using points on a sphere, since for instance $\Kfun(x, y) = \|x-y\|_2^{-1}$ satisfies the potential equation. However, not all kernels satisfy this equation. In particular $\Delta_x \|x-y\|_2^{-2} \neq 0$. As a result, one can't apply Green's theorem, and points on a sphere are not enough to interpolate this kernel. This is illustrated in \autoref{plate-coil computational results}.

A representation similar to \autoref{eq:greenDeri} has been used before.  One idea, the Green hybrid method
(GrH), is explored in \cite{Borm2014GreenHybrid}. This method takes advantage of a two-step strategy. It first analytically approximates the kernel using an integral representation formula, then further compresses its rank by a factor of two using a nested cross approximation. As presented in \cite{Borm20132DSphere}, another way of using Green's formula is to get a parametrization by placing a bounding circle around the geometry, and then spreading the quadrature points across the circle with equal distance.


As shown in \autoref{eq:green}, instead of using both $\Kfun$ and its normal derivative in
\cite{Borm2014GreenHybrid} and \cite{Borm20132DSphere}, our method only uses the kernel
$\Kfun$ itself combined with uniformly distributed points on a bounding sphere. These points
are referred to as ``pseudo-points'' in \cite{Makino1999}.  Similar to our method, \cite{Makino1999}
approximates the potential field using these pseudo-particles to represent the multipole expansion.

In this paper, to deal with 3D geometries, we construct 3D bounding spheres and spread uniformly distributed points as $\Xinit$ and $\Yinit$. The generation of such points is introduced in \autoref{SpherePointsGeneration}.

\subsubsection{Generation} \label{SpherePointsGeneration}

\begin{figure}[!ht]
\centering\includegraphics[width=.3\textwidth]{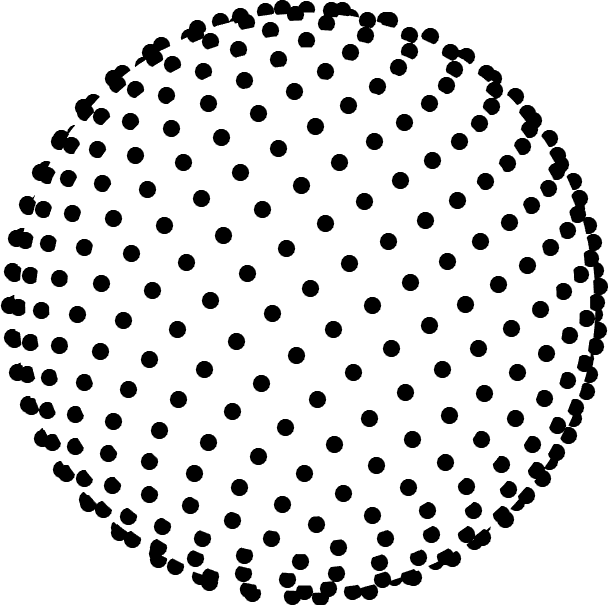}
\caption{Points generated on the surface of a sphere using \autoref{alg:sphere}. Note that the points are uniformly distributed on the spherical surface.}
\label{fig:sphere_gene}
\end{figure}

We use a Fibonacci lattice \cite{gonzalez2010measurement,hardy1979introduction} to place an arbitrary number of points on a sphere around a given cluster of points. See \autoref{fig:sphere_gene} for an illustration of the result. The algorithm \textsf{Sphere}($M$,$n$) constructs a sequence of points $x_i$ using, with $i=1$, \ldots, $n$:
\begin{gather}\label{alg:sphere}
\theta_i = (3-\sqrt{5}) \pi i, \quad
z_i = -1+\frac{2}{n-1} (i-1), \quad
r_i = \sqrt{1-z_i^2} \\
x_i = c_M + r_M
\begin{bmatrix} r_i \cos\theta_i &
r_i \sin\theta_i &
z_i \end{bmatrix}^\top
\end{gather}
where $c_M$ is the center of $M$ and $r_M$ its radius. For this
   case, we use the identity for $\Wx$. \autoref{alg:genericSI} with
                                                     ${\sf
                                                        GenInitSet}(X,r)
                                                        = ({\sf
                                                              Sphere}(X,r),\Wx)$
                                                        is called
                                                        {\sf
                                                           SI-sphere}.

\subsubsection{Discussion}

Compared to the previously introduced method, a major theoretical advantage of distributing points on a sphere around the domain is that they lie on a lower-dimensional manifold. A sphere is only a $d-1$ dimensional manifold in an ambient space of dimension $d$.

\subsection{Random vertices}

Random vertices is the other kind of endo-vertices we used to select $\Xinit$ and
$\Yinit$ directly from $X$ and $Y$. It is a natural and well-known way to
perform sampling from $X$ and $Y$ randomly \cite{halko2011finding,martinsson2016randomized,gu2016efficient,mary2015performance}.

} 

\ifthenelse{1=1}{}{
\item{\bf Further compression the low-rank approximation} 
\par

Given an approximation
\[ \| \KXY - \KXYhat \KXhatYhat \KXhatY \|_F \leq \varepsilon \| \KXY \|_F \]
we point out that one can \emph{always}, no matter the sets $\Xhat$ and $\Yhat$, obtain a re-compressed approximation with a rank equal to or near the optimal $\varepsilon$-rank of $\KXY$.

To see this, consider a matrix $A$ and a matrix $B$ of rank $r_1$ (in our problem, $A = \KXY$ and $B = \KXYhat \KXhatYhat^{-1} \KXhatY$). The key is that
\[ \| A - B \|_2 \leq \varepsilon \Rightarrow | \sigma_k(A) - \sigma_k(B) | \leq \varepsilon. \]
This is a direct consequence of corollary 8.6.2 in \cite{golub13}. Intuitively, this means the singular values are continuous. And this imply that, up to $\varepsilon$, $A$ and $B$ have roughly the same singular values. Hence, we can expect to be able to recompress $B$ to another low-rank matrix with an even smaller rank, $r_2$, closer to the $\varepsilon$-rank of $A$.

In practice, this can be done by computing QR factorizations of the left and right factors, and compressing the central matrix. \autoref{alg:recomp} illustrates the process. We denote by $U_{x|:,1:r_2}$ the matrix obtained by keeping only the first $r_2$ columns of $U_x$, and similarly for $U_y$. $S_{1:r_2, 1:r_2}$ is the leading $r_2 \times r_2$ bloc of $S$.
One can easily verify that
\begin{align*} \KXYhat \KXhatYhat^{-1} \KXhatY & = Q_x R_x \KXhatYhat^{-1} R_y^\top Q_y^\top \\
                                                                        & = Q_x U_x S U_y^\top Q_y^\top  \\
                                                                        & \approx Q_x U_{x|:,1:r_2} S_{1:r_2,1:r_2} U_{y|:,1:r_2} Q_y \end{align*}
where the error in the truncation of $S$ has not been further amplified since $Q_x$ and $Q_y$ are orthogonal.

\begin{algorithm}
\caption{Recompression : $UV^\top = {\sf Recomp}(\Kfun, X, Y, \Xhat, \Yhat, \varepsilon)$}
\label{alg:recomp}
\begin{algorithmic}
\REQUIRE{Approximation $\KXY \approx \KXYhat \KXhatYhat^{-1} \KXhatY$ of rank $r_1$}
\ENSURE{Recompression $\KXY \approx UV^\top$ of rank $r_2 \leq r_1$}
\STATE Compute the QR factorizations
\[ Q_x R_x = \KXYhat, \quad
Q_y R_y = \KXhatY^\top \]

Compute the SVD
\[ U_x S U_y = R_x \KXhatYhat^{-1} R_y^\top \]

Truncate $S$ at $\varepsilon$ and obtain the rank $r_2$
\[ r_2 = \max\{1 \leq i \leq r_1 | S_{11} \varepsilon \leq S_{ii} \} \]

\RETURN $UV^\top$ with $U = Q_x U_{x|:,1:r_2}$ and $V = Q_y U_{y|:,1:r_2} S_{1:r_2,1:r_2}$
\end{algorithmic}
\end{algorithm}

Finally, from the result introduced above, we see that the singular values of $\KXY$ and $U_x S U_y$ are close (up top $\varepsilon$). Hence, by truncating up to $\varepsilon$, we can expect to recover the $\varepsilon$ rank of $\KXY$.
We illustrate this result on a concrete example in \autoref{sec:comp_res}.
\subsection{Computational complexities}

The computational complexities of the Skeletonized Interpolation and the recompression process are summarized below, where we assume $|X|=|Y|=n$. Furthermore, we assume there that the cost of evaluating the kernel on $X \times Y$ is $O(|X| |Y|)$.

For Skeletonized Interpolation, the leading cost is the {\sf GenInitSet} algorithm (if using MDV), the RRQRs and the construction of the left and right factors. Many algorithms exist to compute RRQRs; in our experiments, we use a simple column-pivoted QR algorithm with cost $O(r_0^2 r_1)$ (\cite{golub13}, page 278).


%
%
%


Underline is used to distinguish the number of kernel evaluations (assuming to require $O(1)$ flops) from straight arithmetic operations on floating point numbers.

\medskip

\noindent {\bf Skeletonized Interpolation (\autoref{alg:genericSI})}\\
Maximally-dispersed vertices: $O(n r_0)$; all other methods: $O(r_0)$

Cost of each step:
\begin{itemize}
  \item \underline{Build $K_w$}: $O(r_0^2)$
  \item RRQR of $K_w$ and $K_w^\top$: $O(r_0^2 r_1)$
  \item \underline{Build $\KXYhat$ and $\KXhatY$}: $O(n r_1)$
  \item \underline{Build $\KXhatYhat$}: $O(r_1^2)$
  \item LU factorization of $\KXhatYhat$: $O(r_1^3)$
  \end{itemize}

Total: \text{arithmetic cost:} $O(r_0^2 r_1)$; \text{kernel evaluations:} $O(n r_1 + r_0^2)$

\medskip

\noindent \textbf{Recompression (\autoref{eqn:compress})}
\begin{itemize}
                \item QR factorization of $\KXhatY$ and $\KXYhat$: $O(n r_1^2)$
  \item Compute $R_x \KXhatYhat^{-1} R_y^\top$: $O(r_1^3)$
  \item SVD of $R_x \KXhatYhat^{-1} R_y^\top$: $O(r_1^3)$
  \end{itemize}

Total arithmetic cost: $O(n r_1^2)$
} 
\end{itemize}
%
%

\section{Numerical experiments} \label{Numerical experiments}

\def\HH{6.3cm}
\def\WW{6.8cm}
\def\SS{-2.3ex}
\def\LL{0.7cm}
\def\CC{1.1cm}
We describe numerical experiments on a set of three meshes, a torus,
a plate and coil modeling electromagnetic forming,
and an engine block from acoustic analysis.
We use two kernel functions, $1/r$ and $1/r^2$.
\begin{equation}
\Kfun(x,y) = \frac{1}{\|x-y\|_2}
   \text{\  or\  }
\Kfun(x,y) = \frac{1}{\|x-y\|_2^2}
\end{equation}
The three meshes were generated by LS-DYNA \cite{LSDYNA} from 
Livermore Software Technology Corporation. 
The torus and plate-coil cases are from the electromagnetic solver
(see \cite{EM_BEM2009}) and the engine case is from the acoustics 
solver (see \cite{YunHuangAcousics}).
\begin{figure}[!ht]
\centering
\subfloat[Torus geometry. Inner and outer radii are 3 and 8, respectively\label{fig:torus_geometry}]{\includegraphics[width=0.32\textwidth]{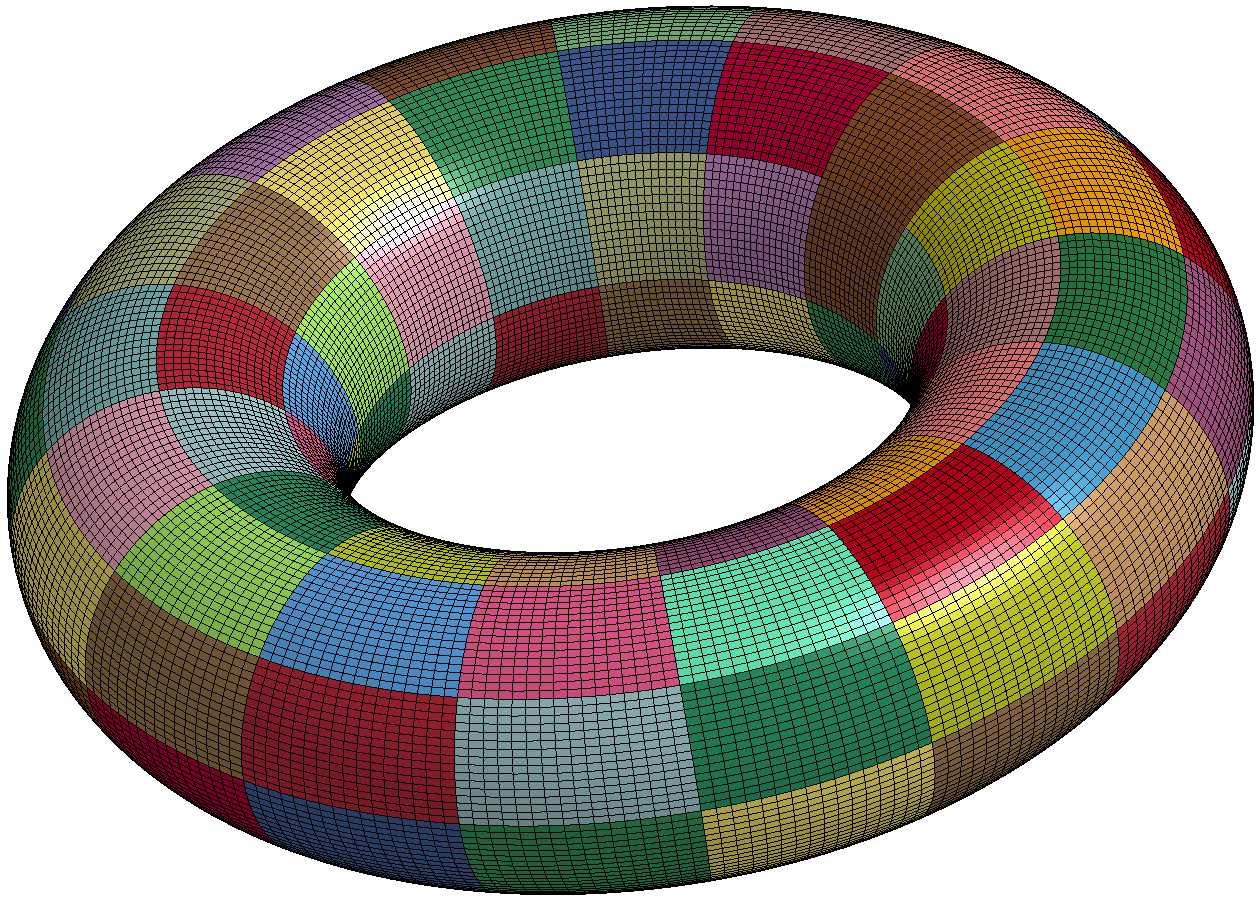}}
\subfloat[Plate-coil geometry\label{fig:coil_geometry}]{\includegraphics[width=0.32\textwidth]{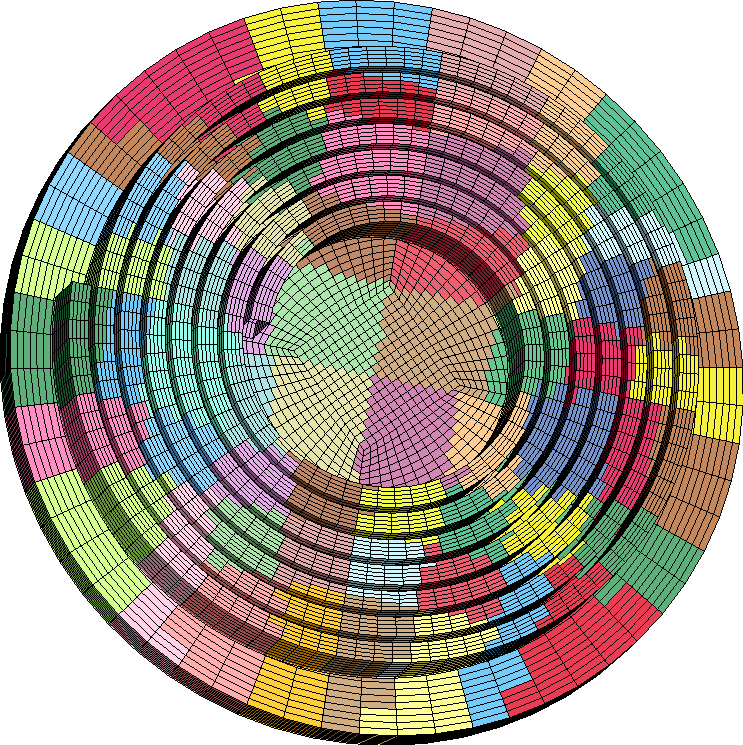}}
\subfloat[Engine geometry\label{fig:engine_geometry}]{\includegraphics[width=0.32\textwidth]{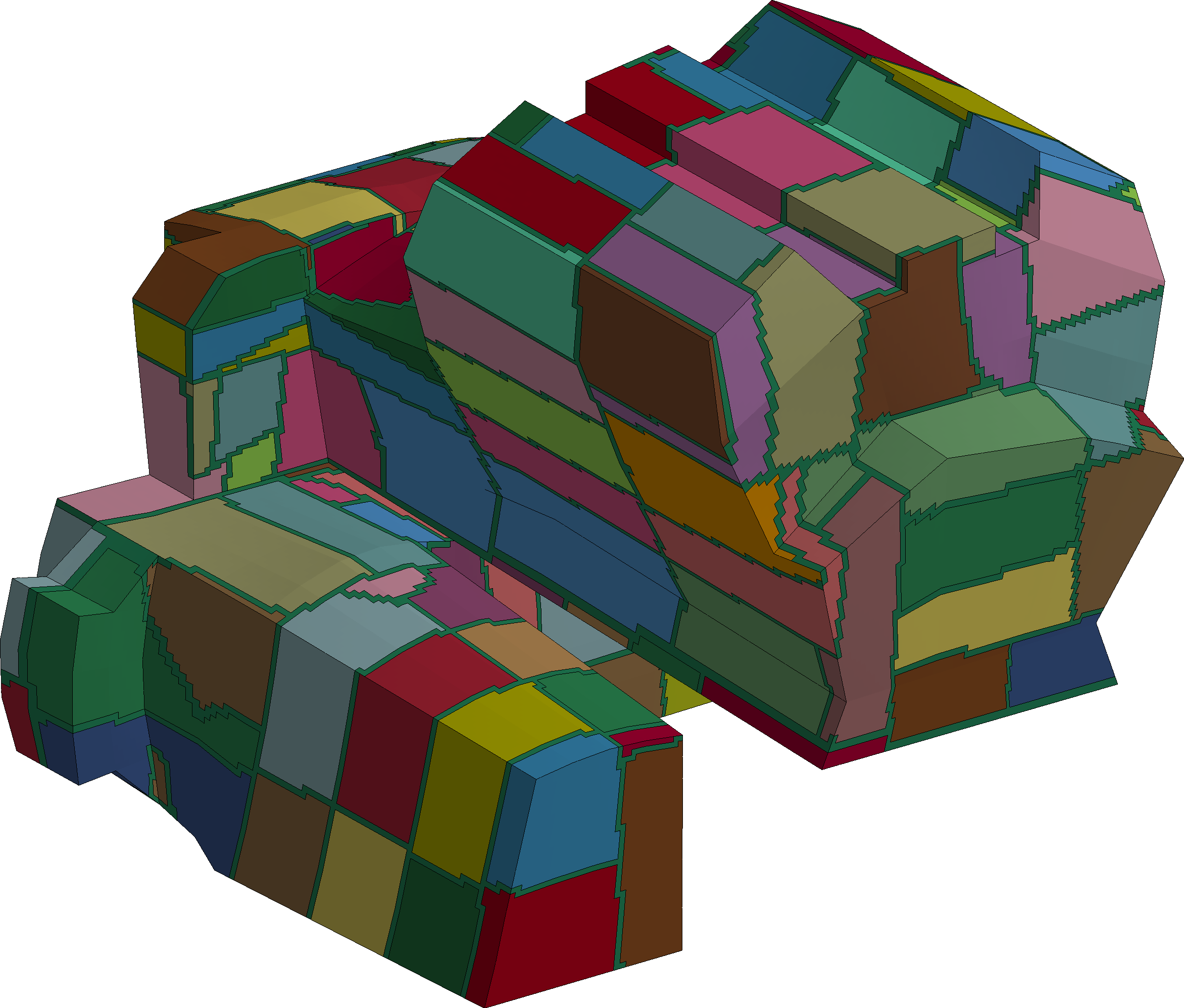}}
\caption{The three benchmarks geometries.}
\label{fig:all_geometries}
\end{figure}
\par
The three meshes are moderate in size, $O(10,000)$ discretization
points, domain size $|X|$, where $250 \le |X| \le 650$.
The standard deviation is a measure of the load imbalance due to
subdomains of different sizes.
Note, the torus mesh has 128 domains of equal size.
\par \smallskip \par
\begin{center}
\begin{tabular}{|lrr|rrrr|}
\hline
\multicolumn{3}{|c|}{} & \multicolumn{4}{c|}{$|X|$} \\
\textbf{Geometry} & \textbf{\# Vertices} & \textbf{\# domains} 
      &  \textbf{Min.} & \textbf{Max.} &  \textbf{Mean}  
      & \textbf{Std Dev.} \\ \hline
{\sf torus}      & 32,768 & 128 & 256 & 256 & 256 &  0.0 \\
{\sf plate-coil} & 20,794 &  64 & 309 & 364 & 325 & 15.4 \\
{\sf engine}     & 78,080 & 128 & 573 & 647 & 610 &  4.6 \\ \hline
\end{tabular}
\end{center}
\par \smallskip \par
We study each of the four algorithms to construct the interpolation
points.
To compute the $\Xinit$, $\Yinit$, $\whX$ and $\whY$ sets,
we increase the sizes of $\Xinit$ and $\Yinit$ by a factor 
$\omega = 1.1$ and perform SI until the relative error 
\begin{equation}
\frac{\left\| \KXY - \KXYhat\KXhatYhat^{-1} \KXhatY 
                        \right\|_F}
                       {\left\|\KXY\right\|_F}
\le \varepsilon^*
\end{equation}
is smaller than a desired tolerance $\varepsilon^*$.
We set the initial $r_0 = 1$ for all strategies but {\sf
   SI-Chebyshev},
for which we use $r_0 = 8$ (2 nodes in each 3 dimensions).
\par
This error tolerance $\varepsilon^*$ is different from the 
$\varepsilon$ in \autoref{alg:genericSI}, which is used to 
extract $\Xhat$ and $\Yhat$ based on the $\varepsilon$-rank.
Here, we use $\varepsilon = 0.1 \varepsilon^*$.
We vary $\varepsilon^* = 10^{-3},\dots,10^{-10}$ and collect 
the resulting ranks $r_0 = |\Xinit| = |\Yinit|$ and $r_1 = |\widehat
X| = |\widehat Y|$.
%
%
\ifthenelse{1=1}{}{
\begin{algorithm}
\caption{Adaptive SI: $\widehat K_{X,Y} = {\sf ASI}(\Kfun, X, Y, \varepsilon^*, \varepsilon, {\sf GenInitSet}, \omega, r_{0})$}
\label{alg:adaptiveSI}
\begin{algorithmic}[1]
\REQUIRE{Kernel $\Kfun$, meshes $X$ and $Y$, approximation tolerance $\varepsilon^*$, RRQR tolerance $\varepsilon$, initial point selection algorithm ${\sf GenInitSet}$, rank growth factor $\omega$, initial rank $r_0$}
\ENSURE{Approximation $\widehat K_{X,Y}$ of $\KXY$}
\WHILE{$\hat\varepsilon > \varepsilon^*$}
    \STATE Using ${\sf SI}$ (\autoref{alg:genericSI}), build the low-rank approximation
    \[ \widehat K_{X,Y} = {\sf SI}(\Kfun, X, Y, \varepsilon, {\sf GenInitSet}, r_0) \]
    and compute the error
    $\hat\varepsilon = \frac{\left\| \widehat K_{X,Y} - \KXY\right\|_F}{\left\|\KXY\right\|_F}$. \\
    \STATE Set $r_0 \leftarrow \omega r_0$.
\ENDWHILE
\RETURN $\widehat K_{X,Y}$
\end{algorithmic}
\end{algorithm}
} 
%
%
We measure the distance between two domains $X_i$ and $X_j$ by
the distance ratio $\dr(i,j)$.
\begin{equation} \label{eqn:dr}
\dr(i,j) 
= \frac{\dist(c_i,c_j)}
       {\min(r_i, r_j)}
= \frac{\dist(\text{centroid}(X_i),\text{centroid}(X_j))}
       {\min(\text{radius}(X_i), \text{radius}(X_j))}
\end{equation}
The four strategies are applied to compute low rank
factorizations for all subdomain pairs $(i,j)$ 
with $\dr(i,j) \geq 1$.
\par
The {\sf torus} mesh has perfect load balance in subdomain size,
while {\sf plate-coil} and {\sf engine} have irregularities,
in size and in shape.
The bounding boxes of {\sf SI-Chebyshev} 
and the spheres of {\sf SI-sphere}
fit closely around the irregular domains.
\par
\par
For the sake of simplicity,
the experiments used the identity matrix as weights matrix for {\sf SI-MDV} and {\sf SI-random}, since they generate roughly uniformly distributed points in the volume.
{\sf SI-Chebyshev} and {\sf SI-sphere} used the weights described in \autoref{Initial Interpolation Points Selection}.
\par
For each mesh and the $1/r$ kernel, we used
{\sf SI-Chebyshev}, {\sf SI-sphere}, {\sf SI-MDV} and {\sf
   SI-random} 
to compute $r_0 = |\Xinit| = |\Yinit|$ and $r_1 = |\whX| = |\whY|$.
We split the submatrix pairs into three sets using the distance
ratio --- near-field, mid-field and far-field submatrices.
For each accuracy, for each set we compute mean values of $r_0$ 
and $r_1$, which we compare to the mean SVD rank at that accuracy.
\par
The results for the $1/r^2$ kernel are very close to that of the
$1/r$ kernel, except for one case.
In \autoref{plate-coil computational results}, we see that 
{\sf SI-sphere} does not work for this kernel, as anticipated.
An enclosing surface does not suffice, 
interior vertices must be included in $\Xinit$.

\ifthenelse{1=1}{}{
\begin{table}[!ht]
\caption{Geometries statistics. Clusters is the number of subdomains. Min, max, mean and standard deviations are relative to the cluster sizes.}
\begin{center}
\begin{tabular}{lllllll}
\toprule
\textbf{Geometry} & \textbf{Vertices} & \textbf{Clusters} &  \textbf{Min.} & \textbf{Max.} &  \textbf{Mean}  & \textbf{Std Dev.} \\
\midrule
Torus  & 32,768      & 128 & 256 & 256 & 256    & 0    \\
Plate-coil &  20,794 &  64 & 309 & 364 & 324.91 & 15.36 \\
Engine & 78,080      & 128 & 573 & 647 & 610    & 4.56 \\
\bottomrule
\end{tabular}
\label{tab:Stats}
\end{center}
\end{table}
} 
\subsection{{\sf torus} : 32,768 vertices, 128 domains}
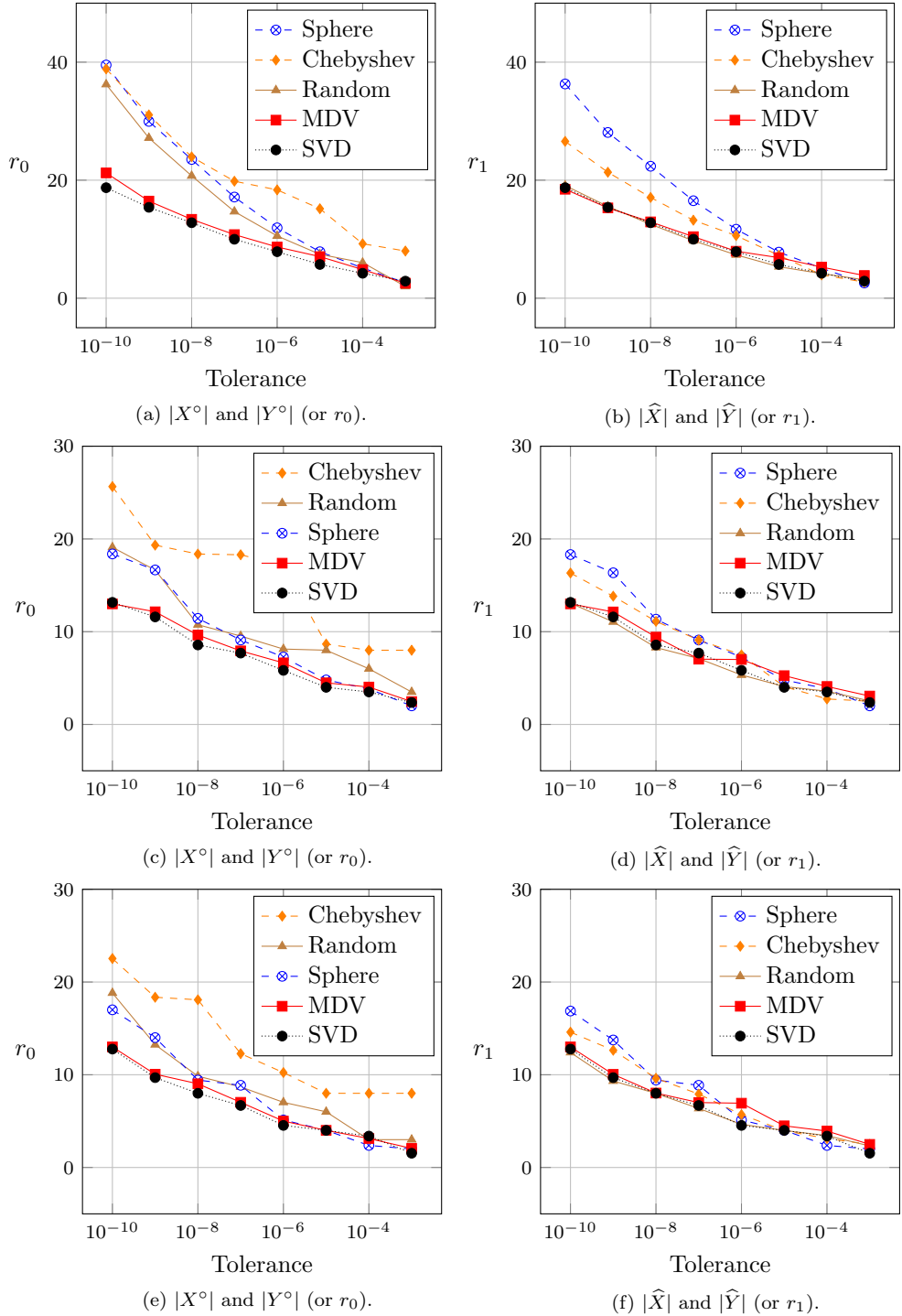
\begin{figure}[]
\hspace{\LL}
\subfloat[$|\Xinit|$ and $|\Yinit|$ (or $r_0$). \label{fig:torus_near_r0}]{
\begin{tikzpicture}[trim axis left,trim axis right]
\begin{semilogxaxis}[
name=plot1,
width=\WW,height=\HH,
ymin = -5, ymax = 50,
xlabel={Tolerance},
ylabel={$r_0$},
grid = major,
legend entries={Sphere, Chebyshev, Random, MDV, SVD},
legend cell align={left},
cycle list name=r0 list near,
ylabel style={rotate=-90},
]
\addplot table [x={Error}, y={r0-sphere}] {plot_torus1_near_r0.dat};
\addplot table [x={Error}, y={r0-Chebyshev}] {plot_torus1_near_r0.dat};
\addplot table [x={Error}, y={r0-random}] {plot_torus1_near_r0.dat};
\addplot table [x={Error}, y={r0-MDV}] {plot_torus1_near_r0.dat};
\addplot table [x={Error}, y={SVD}] {plot_torus1_near_r0.dat};
\end{semilogxaxis}
\end{tikzpicture}
} \hspace{\CC}
\subfloat[$|\Xhat|$ and $|\Yhat|$ (or $r_1$). \label{fig:torus_near_r1}]{
\begin{tikzpicture}[trim axis left,trim axis right]
\begin{semilogxaxis}[
name=plot2,
width=\WW,height=\HH,
at=(plot1.right of south east), anchor=left of south west,
ymin = -5, ymax = 50,
xlabel={Tolerance},
ylabel={$r_1$},
grid = major,
legend entries={Sphere, Chebyshev, Random, MDV, SVD},
legend cell align={left},
cycle list name=r1 list,
ylabel style={rotate=-90},
]
\addplot table [x={Error}, y={r1-sphere}] {plot_torus1_near_r1.dat};
\addplot table [x={Error}, y={r1-Chebyshev}] {plot_torus1_near_r1.dat};
\addplot table [x={Error}, y={r1-random}] {plot_torus1_near_r1.dat};
\addplot table [x={Error}, y={r1-MDV}] {plot_torus1_near_r1.dat};
\addplot table [x={Error}, y={SVD}] {plot_torus1_near_r1.dat};
\end{semilogxaxis}
\end{tikzpicture}
} \\ [\SS] \phantom{} \hspace{\LL}
\subfloat[$|\Xinit|$ and $|\Yinit|$ (or $r_0$). \label{fig:torus_mid_r0}]{
\begin{tikzpicture}[trim axis left,trim axis right]
\begin{semilogxaxis}[
name=plot3,
width=\WW,height=\HH,
ymin = -5, ymax = 30,
xlabel={Tolerance},
ylabel={$r_0$},
grid = major,
legend entries={Chebyshev, Random, Sphere, MDV, SVD},
legend cell align={left},
cycle list name=r0 list mid,
ylabel style={rotate=-90},
]
\addplot table [x={Error}, y={r0-Chebyshev}] {plot_torus1_mid_r0.dat};
\addplot table [x={Error}, y={r0-random}] {plot_torus1_mid_r0.dat};
\addplot table [x={Error}, y={r0-sphere}] {plot_torus1_mid_r0.dat};
\addplot table [x={Error}, y={r0-MDV}] {plot_torus1_mid_r0.dat};
\addplot table [x={Error}, y={SVD}] {plot_torus1_mid_r0.dat};
\end{semilogxaxis}
\end{tikzpicture}
} \hspace{\CC}
\subfloat[$|\Xhat|$ and $|\Yhat|$ (or $r_1$). \label{fig:torus_mid_r1}]{
\begin{tikzpicture}[trim axis left,trim axis right]
\begin{semilogxaxis}[
name=plot4,
width=\WW,height=\HH,
at=(plot3.right of south east), anchor=left of south west,
ymin = -5, ymax = 30,
xlabel={Tolerance},
ylabel={$r_1$},
grid = major,
legend entries={Sphere, Chebyshev, Random, MDV, SVD},
legend cell align={left},
cycle list name=r1 list,
ylabel style={rotate=-90},
]
\addplot table [x={Error}, y={r1-sphere}] {plot_torus1_mid_r1.dat};
\addplot table [x={Error}, y={r1-Chebyshev}] {plot_torus1_mid_r1.dat};
\addplot table [x={Error}, y={r1-random}] {plot_torus1_mid_r1.dat};
\addplot table [x={Error}, y={r1-MDV}] {plot_torus1_mid_r1.dat};
\addplot table [x={Error}, y={SVD}] {plot_torus1_mid_r1.dat};
\end{semilogxaxis}
\end{tikzpicture}
} \\ [\SS]  \phantom{} \hspace{\LL}
\subfloat[$|\Xinit|$ and $|\Yinit|$ (or $r_0$). \label{fig:torus_far_r0}]{
\begin{tikzpicture}[trim axis left,trim axis right]
\begin{semilogxaxis}[
name=plot5,
width=\WW,height=\HH,
ymin = -5, ymax = 30,
xlabel={Tolerance},
ylabel={$r_0$},
grid = major,
legend entries={Chebyshev, Random, Sphere, MDV, SVD},
legend cell align={left},
cycle list name=r0 list mid,
ylabel style={rotate=-90},
]
\addplot table [x={Error}, y={r0-Chebyshev}] {plot_torus1_far_r0.dat};
\addplot table [x={Error}, y={r0-random}] {plot_torus1_far_r0.dat};
\addplot table [x={Error}, y={r0-sphere}] {plot_torus1_far_r0.dat};
\addplot table [x={Error}, y={r0-MDV}] {plot_torus1_far_r0.dat};
\addplot table [x={Error}, y={SVD}] {plot_torus1_far_r0.dat};
\end{semilogxaxis}
\end{tikzpicture}
} \hspace{\CC}
\subfloat[$|\Xhat|$ and $|\Yhat|$ (or $r_1$). \label{fig:torus_far_r1}]{
\begin{tikzpicture}[trim axis left,trim axis right]
\begin{semilogxaxis}[
name=plot6,
width=\WW,height=\HH,
at=(plot5.right of south east), anchor=left of south west,
ymin = -5, ymax = 30,
xlabel={Tolerance},
ylabel={$r_1$},
grid = major,
legend entries={Sphere, Chebyshev, Random, MDV, SVD},
legend cell align={left},
cycle list name=r1 list,
ylabel style={rotate=-90},
]
\addplot table [x={Error}, y={r1-sphere}] {plot_torus1_far_r1.dat};
\addplot table [x={Error}, y={r1-Chebyshev}] {plot_torus1_far_r1.dat};
\addplot table [x={Error}, y={r1-random}] {plot_torus1_far_r1.dat};
\addplot table [x={Error}, y={r1-MDV}] {plot_torus1_far_r1.dat};
\addplot table [x={Error}, y={SVD}] {plot_torus1_far_r1.dat};
\end{semilogxaxis}
\end{tikzpicture}
}
\caption{The mean of the computed ranks in {\sf torus},
$r_0$ on the left, $r_1$ on the right, 
for $2,347$ pairs $\dr \in [1, 4.41]$ (top), 
for $2,347$ pairs $\dr \in [4.41, 6.59]$ (middle) 
and 
for $2,346$ pairs $\dr \in [6.59, 7.87]$ (bottom).}
\label{fig:torus}
\end{figure}

\hfill \break
There are 8128 off-diagonal submatrices in the lower block triangle.
Of these, $7,040$ pairs had distance ratio $\dr > 1$,
which we split into three equally sized sets to represent
near-field, medium-field and far-field pairs.
\Autoref{fig:torus} plots mean values of $r_0$ and $r_1$
for the three different sets of submatrices.
In the torus case, all the
subdomains are of the same shape and all the clusters have a same
distribution. Due to the large distance ratio, SVD ranks are lower
than 20 for all given tolerances.
\pgfplotscreateplotcyclelist{r0 list}{
brown, every mark/.append style={solid, fill=brown}, mark=triangle*\\
dashed, blue, every mark/.append style={solid, fill=white}, mark=otimes*\\
dashed, orange, every mark/.append style={solid, fill=orange},mark=diamond*\\%
red, every mark/.append style={solid, fill=red}, mark=square*\\%
densely dotted, black, every mark/.append style={solid, fill=black}, mark=*\\%
}
\pgfplotscreateplotcyclelist{r0 list near}{
dashed, blue, every mark/.append style={solid, fill=white}, mark=otimes*\\
dashed, orange, every mark/.append style={solid, fill=orange},mark=diamond*\\%
brown, every mark/.append style={solid, fill=brown}, mark=triangle*\\
red, every mark/.append style={solid, fill=red}, mark=square*\\%
densely dotted, black, every mark/.append style={solid, fill=black}, mark=*\\%
}

\pgfplotscreateplotcyclelist{r0 list mid}{
dashed, orange, every mark/.append style={solid, fill=orange},mark=diamond*\\%
brown, every mark/.append style={solid, fill=brown}, mark=triangle*\\
dashed, blue, every mark/.append style={solid, fill=white}, mark=otimes*\\
red, every mark/.append style={solid, fill=red}, mark=square*\\%
densely dotted, black, every mark/.append style={solid, fill=black}, mark=*\\%
}

\pgfplotscreateplotcyclelist{r1 list}{
dashed, blue, every mark/.append style={solid, fill=white}, mark=otimes*\\
dashed, orange, every mark/.append style={solid, fill=orange},mark=diamond*\\%
brown, every mark/.append style={solid, fill=brown}, mark=triangle*\\
red, every mark/.append style={solid, fill=red}, mark=square*\\%
densely dotted, black, every mark/.append style={solid, fill=black}, mark=*\\%
}

We observe that to achieve a tolerance of $\varepsilon^* = 10^{-10}$, we only
require $r_0$ less than $20\%$ of the subdomain size $|X|$. Since the ranks are overall small, $r_0$ grows fairly slowly.

While {\sf SI-Chebyshev} and {\sf SI-sphere} have similar $r_0$ for all tolerances, ``$r_1$ Chebyshev'' is less than ``$r_1$ sphere''
(see \Autoref{fig:torus_near_r1,fig:torus_mid_r1,fig:torus_far_r1}).
This means that using {\sf SI-sphere} is less efficient and {\sf
   SI-Chebyshev} can compress the rank further.

As shown in
\Autoref{fig:torus_near_r0,fig:torus_mid_r0,fig:torus_far_r0},
   ``$r_0$ MDV'' is usually smaller than ``$r_0$ random''.
   Therefore, using {\sf SI-MDV} can substantially reduce the sizes
   of $\Xinit$ and $\Yinit$ thus it is more efficient than {\sf
      SI-random}.
However, \Autoref{fig:torus_mid_r1,fig:torus_far_r1} suggest that sometimes applying RRQR on ``$r_0$ MDV'' is less effective than on ``$r_0$ random''.
When the tolerance is high, ``$r_1$ MDV'' is slightly larger than ``$r_1$ random''. When the tolerance is low,
``$r_0$ MDV'' $\approx$ ``$r_1$ MDV''. This is because of the MDV heuristic
works very well in uniformly distributed clusters thus we already found $\Xhat$ and
$\Yhat$ as the MDV sets and do not need to do RRQR to further compress the rank.

\subsection{{\sf plate-coil} : 20,794 vertices, 64 domains}

\label{plate-coil computational results}
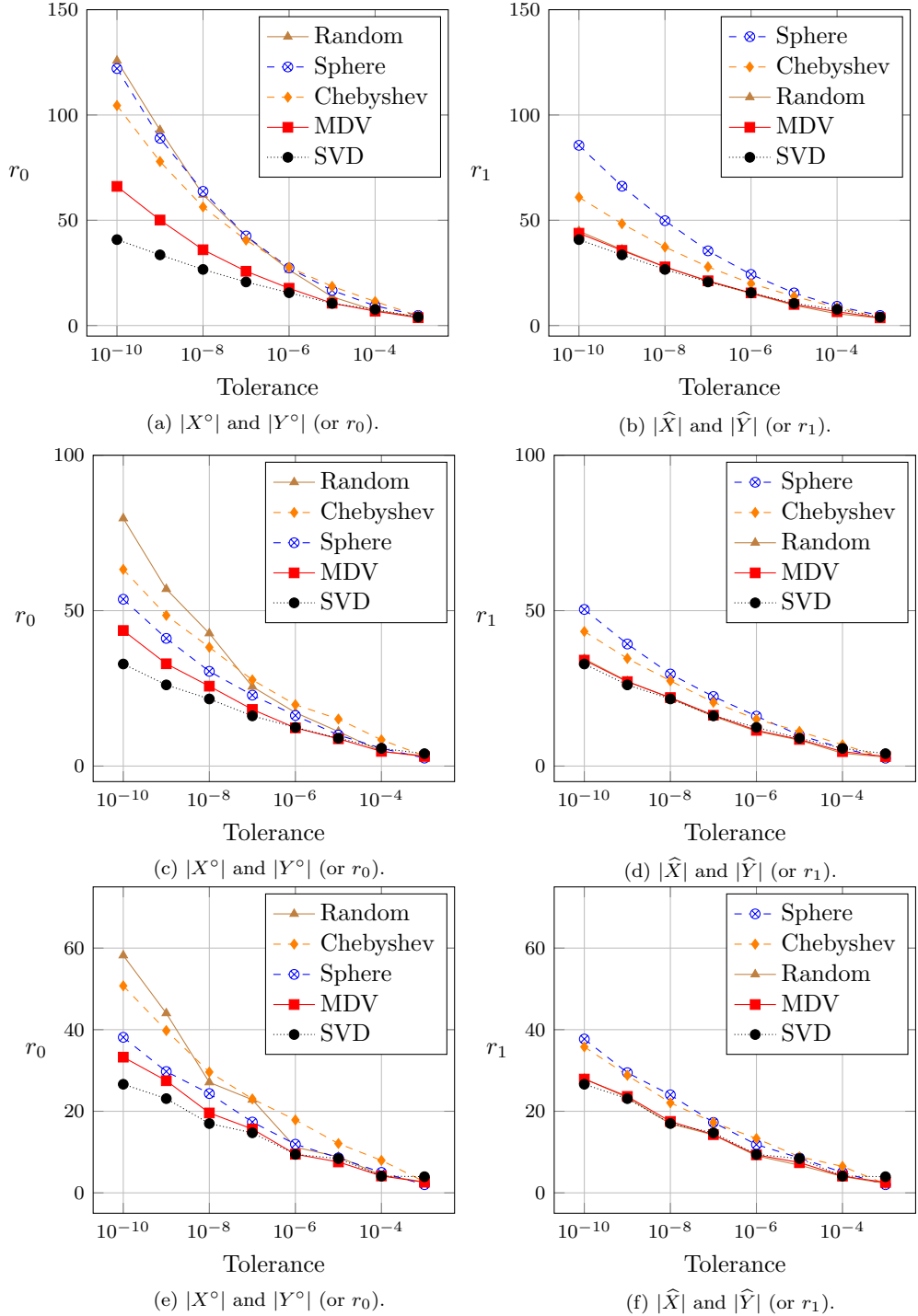
\begin{figure}[]
\hspace{\LL}
\subfloat[$|\Xinit|$ and $|\Yinit|$ (or $r_0$). \label{fig:coil_near_r0}]{
\begin{tikzpicture}[trim axis left,trim axis right]
\begin{semilogxaxis}[
name=plot1,
width=\WW,height=\HH,
ymin = -5, ymax = 150,
xlabel={Tolerance},
ylabel={$r_0$},
grid = major,
legend entries={Random, Sphere, Chebyshev, MDV, SVD},
legend cell align={left},
cycle list name=r0 list,
ylabel style={rotate=-90},
]
\addplot table [x={Error}, y={r0-random}] {plot_coil1_near_r0.dat};
\addplot table [x={Error}, y={r0-sphere}] {plot_coil1_near_r0.dat};
\addplot table [x={Error}, y={r0-Chebyshev}] {plot_coil1_near_r0.dat};
\addplot table [x={Error}, y={r0-MDV}] {plot_coil1_near_r0.dat};
\addplot table [x={Error}, y={SVD}] {plot_coil1_near_r0.dat};
\end{semilogxaxis}
\end{tikzpicture}
} \hspace{\CC}
\subfloat[$|\Xhat|$ and $|\Yhat|$ (or $r_1$). \label{fig:coil_near_r1}]{
\begin{tikzpicture}[trim axis left,trim axis right]
\begin{semilogxaxis}[
name=plot2,
width=\WW,height=\HH,
at=(plot1.right of south east), anchor=left of south west,
ymin = -5, ymax = 150,
xlabel={Tolerance},
ylabel={$r_1$},
grid = major,
legend entries={Sphere, Chebyshev, Random, MDV, SVD},
legend cell align={left},
cycle list name=r1 list,
ylabel style={rotate=-90},
]
\addplot table [x={Error}, y={r1-sphere}] {plot_coil1_near_r1.dat};
\addplot table [x={Error}, y={r1-Chebyshev}] {plot_coil1_near_r1.dat};
\addplot table [x={Error}, y={r1-random}] {plot_coil1_near_r1.dat};
\addplot table [x={Error}, y={r1-MDV}] {plot_coil1_near_r1.dat};
\addplot table [x={Error}, y={SVD}] {plot_coil1_near_r1.dat};
\end{semilogxaxis}
\end{tikzpicture}
}  \\ [\SS] \phantom{} \hspace{\LL}
\subfloat[$|\Xinit|$ and $|\Yinit|$ (or $r_0$). \label{fig:coil_mid_r0}]{
\begin{tikzpicture}[trim axis left,trim axis right]
\begin{semilogxaxis}[
name=plot1,
width=\WW,height=\HH,
ymin = -5, ymax = 100,
xlabel={Tolerance},
ylabel={$r_0$},
grid = major,
legend entries={Random, Chebyshev, Sphere, MDV, SVD},
legend cell align={left},
cycle list name=r0 list coil,
ylabel style={rotate=-90},
]
\addplot table [x={Error}, y={r0-random}] {plot_coil1_mid_r0.dat};
\addplot table [x={Error}, y={r0-Chebyshev}] {plot_coil1_mid_r0.dat};
\addplot table [x={Error}, y={r0-sphere}] {plot_coil1_mid_r0.dat};
\addplot table [x={Error}, y={r0-MDV}] {plot_coil1_mid_r0.dat};
\addplot table [x={Error}, y={SVD}] {plot_coil1_mid_r0.dat};
\end{semilogxaxis}
\end{tikzpicture}
} \hspace{\CC}
\subfloat[$|\Xhat|$ and $|\Yhat|$ (or $r_1$). \label{fig:coil_mid_r1}]{
\begin{tikzpicture}[trim axis left,trim axis right]
\begin{semilogxaxis}[
name=plot2,
width=\WW,height=\HH,
at=(plot1.right of south east), anchor=left of south west,
ymin = -5, ymax = 100,
xlabel={Tolerance},
ylabel={$r_1$},
grid = major,
legend entries={Sphere, Chebyshev, Random, MDV, SVD},
legend cell align={left},
cycle list name=r1 list,
ylabel style={rotate=-90},
]
\addplot table [x={Error}, y={r1-sphere}] {plot_coil1_mid_r1.dat};
\addplot table [x={Error}, y={r1-Chebyshev}] {plot_coil1_mid_r1.dat};
\addplot table [x={Error}, y={r1-random}] {plot_coil1_mid_r1.dat};
\addplot table [x={Error}, y={r1-MDV}] {plot_coil1_mid_r1.dat};
\addplot table [x={Error}, y={SVD}] {plot_coil1_mid_r1.dat};
\end{semilogxaxis}
\end{tikzpicture}
} \\ [\SS] \phantom{} \hspace{\LL}
\subfloat[$|\Xinit|$ and $|\Yinit|$ (or $r_0$). \label{fig:coil_far_r0}]{
\begin{tikzpicture}[trim axis left,trim axis right]
\begin{semilogxaxis}[
name=plot1,
width=\WW,height=\HH,
ymin = -5, ymax = 75,
xlabel={Tolerance},
ylabel={$r_0$},
grid = major,
legend entries={Random, Chebyshev, Sphere, MDV, SVD},
legend cell align={left},
cycle list name=r0 list coil,
ylabel style={rotate=-90},
]
\addplot table [x={Error}, y={r0-random}] {plot_coil1_far_r0.dat};
\addplot table [x={Error}, y={r0-Chebyshev}] {plot_coil1_far_r0.dat};
\addplot table [x={Error}, y={r0-sphere}] {plot_coil1_far_r0.dat};
\addplot table [x={Error}, y={r0-MDV}] {plot_coil1_far_r0.dat};
\addplot table [x={Error}, y={SVD}] {plot_coil1_far_r0.dat};
\end{semilogxaxis}
\end{tikzpicture}
} \hspace{\CC}
\subfloat[$|\Xhat|$ and $|\Yhat|$ (or $r_1$). \label{fig:coil_far_r1}]{
\begin{tikzpicture}[trim axis left,trim axis right]
\begin{semilogxaxis}[
name=plot2,
width=\WW,height=\HH,
at=(plot1.right of south east), anchor=left of south west,
ymin = -5, ymax = 75,
xlabel={Tolerance},
ylabel={$r_1$},
grid = major,
legend entries={Sphere, Chebyshev, Random, MDV, SVD},
legend cell align={left},
cycle list name=r1 list,
ylabel style={rotate=-90},
]
\addplot table [x={Error}, y={r1-sphere}] {plot_coil1_far_r1.dat};
\addplot table [x={Error}, y={r1-Chebyshev}] {plot_coil1_far_r1.dat};
\addplot table [x={Error}, y={r1-random}] {plot_coil1_far_r1.dat};
\addplot table [x={Error}, y={r1-MDV}] {plot_coil1_far_r1.dat};
\addplot table [x={Error}, y={SVD}] {plot_coil1_far_r1.dat};
\end{semilogxaxis}
\end{tikzpicture}
}
\caption{The mean of the computed ranks for $404$ pairs in the
   plate-coil case, $dr \in [1, 1.58]$ (top), $\dr \in [1.58, 2.21]$
      (middle) and $\dr \in [2.21, 3.28]$ (top).
      In~\autoref{fig:coil_near_r1} and \autoref{fig:coil_far_r1},
      {\sf SI-random} is behind {\sf SI-MDV}.}
\label{fig:plate-coil}
\end{figure}

\hfill \break
There are 2,016 off-diagonal submatrices, of these 1,212 pairs
have distance ratio greater than one.
\Autoref{fig:plate-coil} shows plots of mean values of $r_0$ and $r_1$ 
for near-field, medium-field and far-field submatrices.
\par
Unlike {\sf torus}, the subdomains in {\sf plate-coil} 
are irregularly shaped, and the domains do not have uniformly 
distributed points. 
Because {\sf plate-coil} mesh is smaller in size and 
in the number of domains, the distances are shorter.
The closer two subdomains, the higher the SVD rank.
This is also reflected in the size of the initial sets 
$r_0 = |\Xinit|$, up to half the size of the domain for
high accuracy near-field submatrices.

As the distance ratio increases, 
both $r_0$ and $r_1$ decrease for all four strategies.  
The relative performance of the four
strategies is similar to that of torus case.  
Although in
\Autoref{fig:coil_mid_r0,fig:coil_far_r0} ``$r_0$ Chebyshev'' is
higher than ``$r_0$ Sphere'', ``$r_1$ Chebyshev'' is always less
than ``$r_1$ Sphere'' (see
      \Autoref{fig:coil_near_r1,fig:coil_mid_r1,fig:coil_far_r1}).
Thus {\sf SI-Chebyshev} can compress the rank further than {\sf
   SI-sphere}.
\par
As for {\sf SI-MDV} and {\sf SI-random}
(\Autoref{fig:coil_near_r0,fig:coil_mid_r0,fig:coil_far_r0}),
``$r_0$ MDV'' is much smaller than ``$r_0$ random'' thus {\sf
   {\sf SI-MDV}}
is more efficient. 
The final ranks of {\sf SI-MDV} and {\sf SI-random} are nearly equal
to the svd ranks, and further compression is not needed.
However, the final ranks of {\sf SI-Chebyshev} and {\sf SI-sphere} 
are greater than the svd ranks, and further compression is needed.
\par


\pgfplotscreateplotcyclelist{r0 list}{
brown, every mark/.append style={solid, fill=brown}, mark=triangle*\\
dashed, blue, every mark/.append style={solid, fill=white}, mark=otimes*\\
dashed, orange, every mark/.append style={solid, fill=orange},mark=diamond*\\%
red, every mark/.append style={solid, fill=red}, mark=square*\\%
densely dotted, black, every mark/.append style={solid, fill=black}, mark=*\\%
}

\pgfplotscreateplotcyclelist{r0 list coil}{
brown, every mark/.append style={solid, fill=brown}, mark=triangle*\\
dashed, orange, every mark/.append style={solid, fill=orange},mark=diamond*\\%
dashed, blue, every mark/.append style={solid, fill=white}, mark=otimes*\\
red, every mark/.append style={solid, fill=red}, mark=square*\\%
densely dotted, black, every mark/.append style={solid, fill=black}, mark=*\\%
}


\begin{figure}[!ht]
\centering
\begin{tikzpicture}
\begin{semilogxaxis}[
xlabel={Tolerance},
ylabel={$r_0$},
ymin=-5, ymax = 120,
width=9cm,
height=6cm,
grid = major,
legend entries={Sphere, Chebyshev, Random, MDV, SVD},
legend pos = outer north east,
legend cell align={left},
cycle list name=r0 list 1/r2,
ylabel style={rotate=-90},
]
\addplot table [x={Error}, y={r0-sphere}] {plot_coil2_far_r0.dat};
\addplot table [x={Error}, y={r0-Chebyshev}] {plot_coil2_far_r0.dat};
\addplot table [x={Error}, y={r0-random}] {plot_coil2_far_r0.dat};
\addplot table [x={Error}, y={r0-MDV}] {plot_coil2_far_r0.dat};
\addplot table [x={Error}, y={SVD}] {plot_coil2_far_r0.dat};
\end{semilogxaxis}
\end{tikzpicture}
\caption{The mean of computed rank $r_0$ for $404$ far-field pairs
   in plate-coil case, $\dr \in [2.21, 3.28]$. The kernel is
      $\Kfun(x,y) = \|x-y\|_2^{-2}$. Notice for this kernel, {\sf
         SI-sphere} does not work and $r_0$ for {\sf SI-sphere} blows up.}
\label{fig:kernel2}
\end{figure}
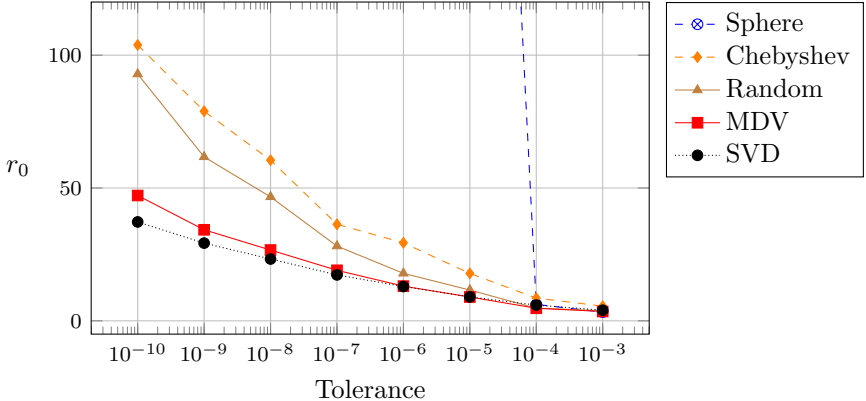

\autoref{fig:kernel2} shows the computational results with the kernel $\Kfun(x,y) = \|x-y\|_2^{-2}$.
The curve of {\sf SI-sphere} is only partially shown. This is because when the given tolerance is less than
$10^{-4}$, $r_0$ of {\sf SI-sphere} blows up and exceeds 500 in our experiments. This means that
building $\Xinit$ and $\Yinit$ using vertices on a sphere does not work with this kernel, as explained in \autoref{sec:sphere}.

\pgfplotscreateplotcyclelist{r0 list 1/r2}{
dashed, blue, every mark/.append style={solid, fill=white}, mark=otimes*\\
dashed, orange, every mark/.append style={solid, fill=orange},mark=diamond*\\%
brown, every mark/.append style={solid, fill=brown}, mark=triangle*\\
red, every mark/.append style={solid, fill=red}, mark=square*\\%
densely dotted, black, every mark/.append style={solid, fill=black}, mark=*\\%
}

\subsection{{\sf engine} : 78,080 vertices, 128 domains}

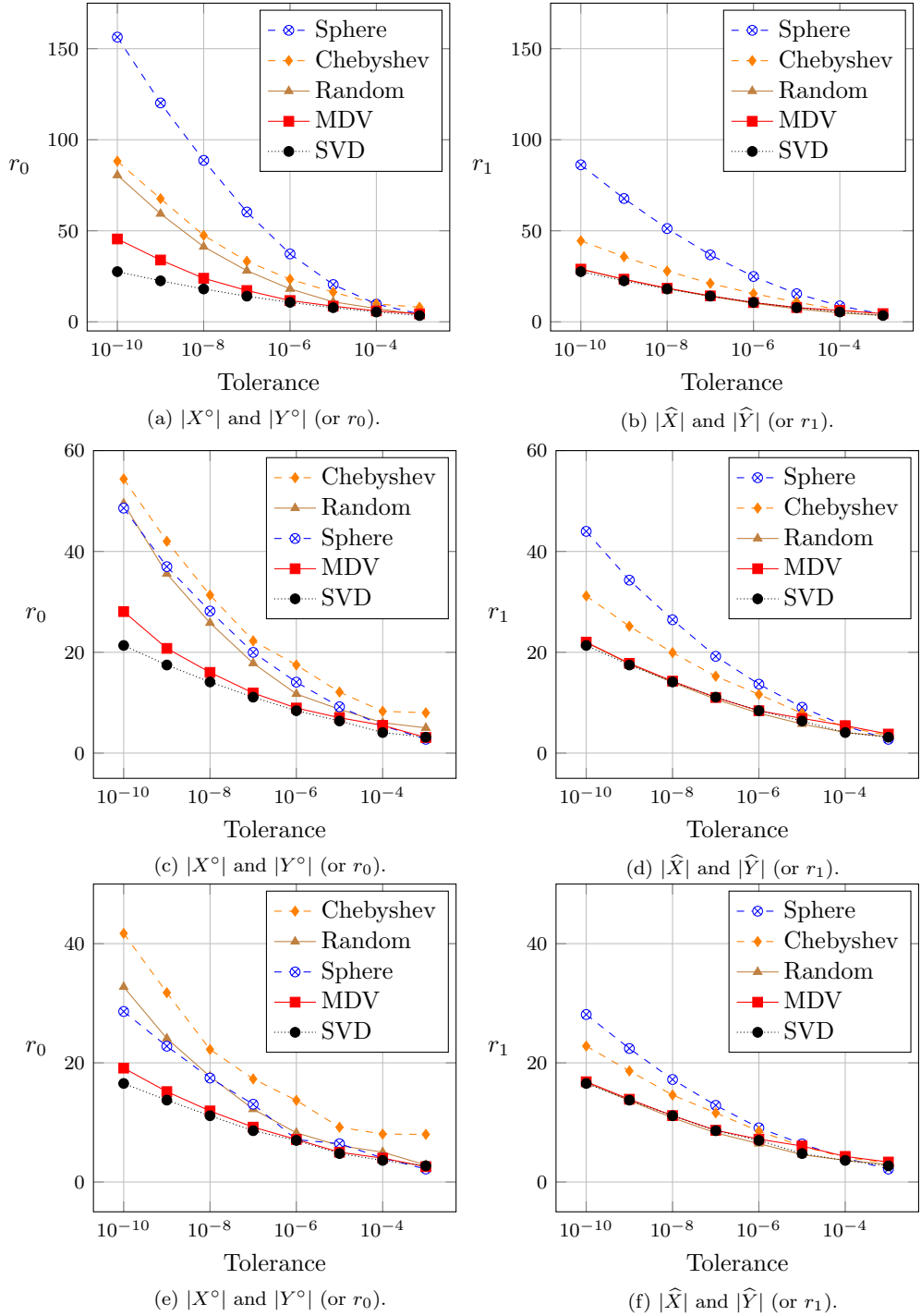
\begin{figure}[]
\hspace{\LL}
\subfloat[$|\Xinit|$ and $|\Yinit|$ (or $r_0$). \label{fig:engine_near_r0}]{
\begin{tikzpicture}[trim axis left,trim axis right]
\begin{semilogxaxis}[
name=plot1,
width=\WW,height=\HH,
ymin = -5, ymax = 175,
xlabel={Tolerance},
ylabel={$r_0$},
grid = major,
legend entries={Sphere, Chebyshev, Random, MDV, SVD},
legend cell align={left},
cycle list name=r0 list near,
ylabel style={rotate=-90},
]
\addplot table [x={Error}, y={r0-sphere}] {plot_engine1_near_r0.dat};
\addplot table [x={Error}, y={r0-Chebyshev}] {plot_engine1_near_r0.dat};
\addplot table [x={Error}, y={r0-random}] {plot_engine1_near_r0.dat};
\addplot table [x={Error}, y={r0-MDV}] {plot_engine1_near_r0.dat};
\addplot table [x={Error}, y={SVD}] {plot_engine1_near_r0.dat};
\end{semilogxaxis}
\end{tikzpicture}
} \hspace{\CC}
\subfloat[$|\Xhat|$ and $|\Yhat|$ (or $r_1$). \label{fig:engine_near_r1}]{
\begin{tikzpicture}[trim axis left,trim axis right]
\begin{semilogxaxis}[
name=plot2,
width=\WW,height=\HH,
at=(plot1.right of south east), anchor=left of south west,
ymin = -5, ymax = 175,
xlabel={Tolerance},
ylabel={$r_1$},
grid = major,
legend entries={Sphere, Chebyshev, Random, MDV, SVD},
legend cell align={left},
cycle list name=r1 list,
ylabel style={rotate=-90},
]
\addplot table [x={Error}, y={r1-sphere}] {plot_engine1_near_r1.dat};
\addplot table [x={Error}, y={r1-Chebyshev}] {plot_engine1_near_r1.dat};
\addplot table [x={Error}, y={r1-random}] {plot_engine1_near_r1.dat};
\addplot table [x={Error}, y={r1-MDV}] {plot_engine1_near_r1.dat};
\addplot table [x={Error}, y={SVD}] {plot_engine1_near_r1.dat};
\end{semilogxaxis}
\end{tikzpicture}
} \\ [\SS] \phantom{} \hspace{\LL}
\subfloat[$|\Xinit|$ and $|\Yinit|$ (or $r_0$). \label{fig:engine_mid_r0}]{
\begin{tikzpicture}[trim axis left,trim axis right]
\begin{semilogxaxis}[
name=plot3,
width=\WW,height=\HH,
ymin = -5, ymax = 60,
xlabel={Tolerance},
ylabel={$r_0$},
grid = major,
legend entries={Chebyshev, Random, Sphere, MDV, SVD},
legend cell align={left},
cycle list name=r0 list mid,
ylabel style={rotate=-90},
]
\addplot table [x={Error}, y={r0-Chebyshev}] {plot_engine1_mid_r0.dat};
\addplot table [x={Error}, y={r0-random}] {plot_engine1_mid_r0.dat};
\addplot table [x={Error}, y={r0-sphere}] {plot_engine1_mid_r0.dat};
\addplot table [x={Error}, y={r0-MDV}] {plot_engine1_mid_r0.dat};
\addplot table [x={Error}, y={SVD}] {plot_engine1_mid_r0.dat};
\end{semilogxaxis}
\end{tikzpicture}
} \hspace{\CC}
\subfloat[$|\Xhat|$ and $|\Yhat|$ (or $r_1$). \label{fig:engine_mid_r1}]{
\begin{tikzpicture}[trim axis left,trim axis right]
\begin{semilogxaxis}[
name=plot4,
width=\WW,height=\HH,
at=(plot3.right of south east), anchor=left of south west,
ymin = -5, ymax = 60,
xlabel={Tolerance},
ylabel={$r_1$},
grid = major,
legend entries={Sphere, Chebyshev, Random, MDV, SVD},
legend cell align={left},
cycle list name=r1 list,
ylabel style={rotate=-90},
]
\addplot table [x={Error}, y={r1-sphere}] {plot_engine1_mid_r1.dat};
\addplot table [x={Error}, y={r1-Chebyshev}] {plot_engine1_mid_r1.dat};
\addplot table [x={Error}, y={r1-random}] {plot_engine1_mid_r1.dat};
\addplot table [x={Error}, y={r1-MDV}] {plot_engine1_mid_r1.dat};
\addplot table [x={Error}, y={SVD}] {plot_engine1_mid_r1.dat};
\end{semilogxaxis}
\end{tikzpicture}
} \\ [\SS] \phantom{} \hspace{\LL}
\subfloat[$|\Xinit|$ and $|\Yinit|$ (or $r_0$). \label{fig:engine_far_r0}]{
\begin{tikzpicture}[trim axis left,trim axis right]
\begin{semilogxaxis}[
name=plot5,
width=\WW,height=\HH,
ymin = -5, ymax = 50,
xlabel={Tolerance},
ylabel={$r_0$},
grid = major,
legend entries={Chebyshev, Random, Sphere, MDV, SVD},
legend cell align={left},
cycle list name=r0 list mid,
ylabel style={rotate=-90},
]
\addplot table [x={Error}, y={r0-Chebyshev}] {plot_engine1_far_r0.dat};
\addplot table [x={Error}, y={r0-random}] {plot_engine1_far_r0.dat};
\addplot table [x={Error}, y={r0-sphere}] {plot_engine1_far_r0.dat};
\addplot table [x={Error}, y={r0-MDV}] {plot_engine1_far_r0.dat};
\addplot table [x={Error}, y={SVD}] {plot_engine1_far_r0.dat};
\end{semilogxaxis}
\end{tikzpicture}
} \hspace{\CC}
\subfloat[$|\Xhat|$ and $|\Yhat|$ (or $r_1$). \label{fig:engine_far_r1}]{
\begin{tikzpicture}[trim axis left,trim axis right]
\begin{semilogxaxis}[
name=plot6,
width=\WW,height=\HH,
at=(plot5.right of south east), anchor=left of south west,
ymin = -5, ymax = 50,
xlabel={Tolerance},
ylabel={$r_1$},
grid = major,
legend entries={Sphere, Chebyshev, Random, MDV, SVD},
legend cell align={left},
cycle list name=r1 list,
ylabel style={rotate=-90},
]
\addplot table [x={Error}, y={r1-sphere}] {plot_engine1_far_r1.dat};
\addplot table [x={Error}, y={r1-Chebyshev}] {plot_engine1_far_r1.dat};
\addplot table [x={Error}, y={r1-random}] {plot_engine1_far_r1.dat};
\addplot table [x={Error}, y={r1-MDV}] {plot_engine1_far_r1.dat};
\addplot table [x={Error}, y={SVD}] {plot_engine1_far_r1.dat};
\end{semilogxaxis}
\end{tikzpicture}
}
\caption{The mean of the computed ranks for $2150$ near-field pairs
   in the engine case, $\dr \in [1, 1.87]$ (top), $\dr \in [1.87,
      2.84]$ (middle) and $\dr \in [2.84, 7.46]$ (bottom). In
         \autoref{fig:engine_near_r1}, {\sf SI-random} is hidden behind
         {\sf SI-MDV} and SVD.}
\label{fig:engine}
\end{figure}

The results are depicted on \autoref{fig:engine}.
As for {\sf plate-coil}, the domains are irregular in size and shape.
However, $r_0$ here is less than that of the plate-coil case 
(about $20\%$ of the subdomain size when the tolerance is low), 
due to the larger distance ratios. 
The relative performance of the four strategies is very similar 
to {\sf plate-coil}. 

\pgfplotscreateplotcyclelist{r0 list near}{
dashed, blue, every mark/.append style={solid, fill=white}, mark=otimes*\\
dashed, orange, every mark/.append style={solid, fill=orange},mark=diamond*\\%
brown, every mark/.append style={solid, fill=brown}, mark=triangle*\\
red, every mark/.append style={solid, fill=red}, mark=square*\\%
densely dotted, black, every mark/.append style={solid, fill=black}, mark=*\\%
}

\pgfplotscreateplotcyclelist{r0 list mid}{
dashed, orange, every mark/.append style={solid, fill=orange},mark=diamond*\\%
brown, every mark/.append style={solid, fill=brown}, mark=triangle*\\
dashed, blue, every mark/.append style={solid, fill=white}, mark=otimes*\\
red, every mark/.append style={solid, fill=red}, mark=square*\\%
densely dotted, black, every mark/.append style={solid, fill=black}, mark=*\\%
}


\subsection{Cases where MDV is inefficient}

The drawback of MDV is that it \textit{a priori} selects points
uniformly in the volume (or using a graph-distance). In most cases,
this ends up working very well. However, for kernels that
are also Green's functions (e.g., when solving some
integral equations), it may be sufficient to place
interpolation points on an enclosing surface. Consider for
example $\log(\|x-y\|)$ in 2D.

We know that this kernel can be approximated using multipole
functions. To obtain a multipole expansion of order $p$, we need
only $O(p)$ coefficients~\cite{Makino1999}; this is in contrast to
general expansions such as Taylor or Chebyshev that require $O(p^2)$
terms. Moreover, one can always find a set of equivalent $p$
``charges'' on an enclosing circle that lead to the same first $p$
terms in the multipole expansion~\cite{Makino1999}. This implies
that to reach order $p$, {\sf SI-sphere} requires $O(p)$ points
only.

\begin{figure}[!ht]
\subfloat[$|\Xinit|$ and $|\Yinit|$ (or $r_0$). \label{fig:mdv_fails_r0}]{
\begin{tikzpicture}
\begin{semilogxaxis}[
name=plot1,
width=6.5cm,height=6.5cm,
ymin = -5, ymax = 80,
xlabel={Tolerance},
ylabel={$r_0$},
grid = major,
legend entries={MDV, Sphere, SVD},
legend cell align={left},
cycle list name=r0 list,
ylabel style={rotate=-90},
]
\addplot table [x={Error}, y={r0-mdv}] {r0_mdv_vs_circle.dat};
\addplot table [x={Error}, y={r0-sphere}] {r0_mdv_vs_circle.dat};
\addplot table [x={Error}, y={SVD}] {r0_mdv_vs_circle.dat};
\end{semilogxaxis}
\end{tikzpicture}
}
\subfloat[$|\Xhat|$ and $|\Yhat|$ (or $r_1$). \label{fig:mdv_fails_r1}]{
\begin{tikzpicture}
\begin{semilogxaxis}[
name=plot1,
width=6.5cm,height=6.5cm,
ymin = -5, ymax = 80,
xlabel={Tolerance},
ylabel={$r_1$},
grid = major,
legend entries={MDV, Sphere, SVD},
legend cell align={left},
cycle list name=r0 list,
ylabel style={rotate=-90},
]
\addplot table [x={Error}, y={r0-mdv}] {r1_mdv_vs_circle.dat};
\addplot table [x={Error}, y={r0-sphere}] {r1_mdv_vs_circle.dat};
\addplot table [x={Error}, y={SVD}] {r1_mdv_vs_circle.dat};
\end{semilogxaxis}
\end{tikzpicture}
} \\
\centering
\subfloat[The geometry and the initial $\Xinit$ and $\Yinit$ for $\varepsilon = 10^{-14}$.]{
\begin{tikzpicture}
\begin{axis}[
name=plotGeo,
width=7cm,height=7cm,
ymin = -1.5, ymax = 4,
xmin = -1.5, xmax = 4,
ytick=\empty,
xtick=\empty,
legend entries={,,$\Xinit_{\text{Sphere}}$, $\Yinit_{\text{Sphere}}$, $\Xinit_{\text{MDV}}$, $\Yinit_{\text{MDV}}$},
legend cell align={left},
legend style={at={(0.325,0.979)}},
ylabel style={rotate=-90},
]
\addplot [black, dashed] table [x={x1}, y={x2}] {geo_mdv_vs_circle.dat};
\addplot [black, dashed] table [x={y1}, y={y2}] {geo_mdv_vs_circle.dat};
\addplot [blue, only marks, mark=o]   table [x={x1}, y={x2}] {circle_mdv_vs_circle_0.dat};
\addplot [blue, only marks, mark=o]   table [x={y1}, y={y2}] {circle_mdv_vs_circle_0.dat};
\addplot [red, only marks, mark=diamond]    table [x={x1}, y={x2}] {mdv_mdv_vs_circle_0.dat};
\addplot [red, only marks, mark=diamond]    table [x={y1}, y={y2}] {mdv_mdv_vs_circle_0.dat};
\end{axis}
\end{tikzpicture}
}  \quad
\subfloat[The geometry and the outputted $\Xhat$ and $\Yhat$ for $\varepsilon = 10^{-14}$. We observe that the points selected by MDV (diamonds) also tend to cluster at the boundaries.\label{fig:mdv_fails_geo}]{
\begin{tikzpicture}
\begin{axis}[
name=plotGeo,
width=7cm,height=7cm,
ymin = -1.5, ymax = 4,
xmin = -1.5, xmax = 4,
ytick=\empty,
xtick=\empty,
legend entries={,,$\Xhat_{\text{Sphere}}$, $\Yhat_{\text{Sphere}}$, $\Xhat_{\text{MDV}}$, $\Yhat_{\text{MDV}}$},
legend cell align={left},
legend style={at={(0.325,0.979)}},
ylabel style={rotate=-90},
]
\addplot [black, dashed] table [x={x1}, y={x2}] {geo_mdv_vs_circle.dat};
\addplot [black, dashed] table [x={y1}, y={y2}] {geo_mdv_vs_circle.dat};
\addplot [blue, only marks, mark=otimes*]   table [x={x1}, y={x2}] {circle_mdv_vs_circle.dat};
\addplot [blue, only marks, mark=otimes*]   table [x={y1}, y={y2}] {circle_mdv_vs_circle.dat};
\addplot [red, only marks, mark=diamond*]    table [x={x1}, y={x2}] {mdv_mdv_vs_circle.dat};
\addplot [red, only marks, mark=diamond*]    table [x={y1}, y={y2}] {mdv_mdv_vs_circle.dat};

\end{axis}
\end{tikzpicture}
} \\

\caption{\label{fig:mdv_fail} Example where MDV does worse than the
points on a sphere. In this geometry, we arrange 833 points
uniformly in two 2D disks of the same shape, separated by less
than one diameter. The kernel is the 2D Laplacian kernel,
$\Kfun(x, y) = \log(\|x-y\|_2)$. We observe than
arranging points on a circle leads to smaller $\Xinit$
and $\Yinit$ for a given accuracy. Chebyshev and
Random give similar results (in terms of ranks $r_0$)
than  MDV points and vertices on a sphere,
respectively, and are omitted for clarity.}
\end{figure}
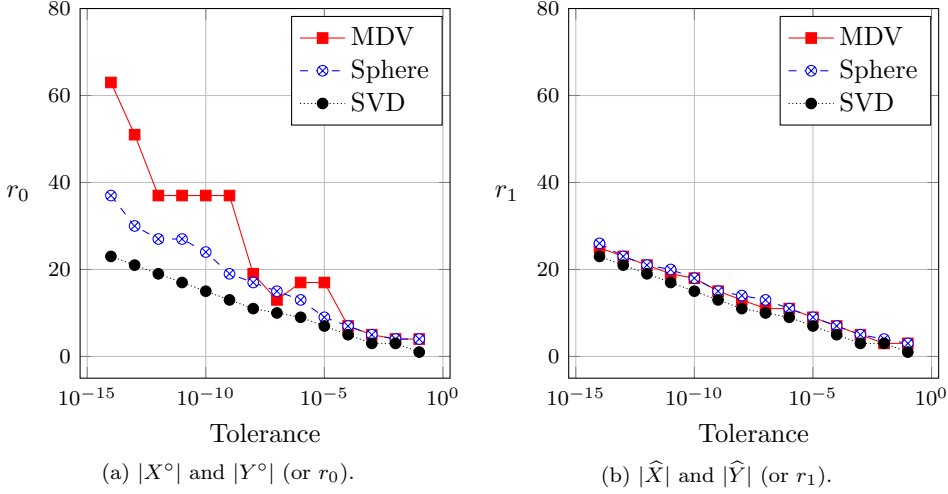
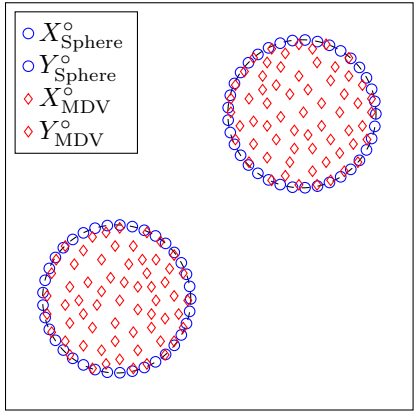
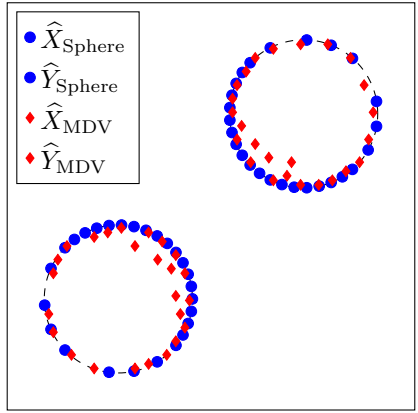

{\sf SI-MDV} requires about the same density of points. The problem
is that if the points in $X$ and $Y$ are uniformly distributed in
the volume, then in order to get the correct density of points
\textbf{near the boundary} one needs $r_0 = O(p^2)$ points in {\sf
SI-MDV}. This is reflected by a much larger value of $r_0$ in the
case of {\sf SI-MDV} compared to {\sf SI-sphere}. See
\autoref{fig:mdv_fail} for an illustrative benchmark. In
particular, see the distribution of points for MDV after RRQR
($\Xhat$ and $\Yhat$), and how it tries to approximate the
distribution from Sphere.

\pgfplotscreateplotcyclelist{r0 list}{
red, every mark/.append style={solid, fill=red}, mark=square*\\%
dashed, blue, every mark/.append style={solid, fill=white}, mark=otimes*\\
densely dotted, black, every mark/.append style={solid, fill=black}, mark=*\\%
}

\subsection{Example of recompression} \label{sec:comp_res}
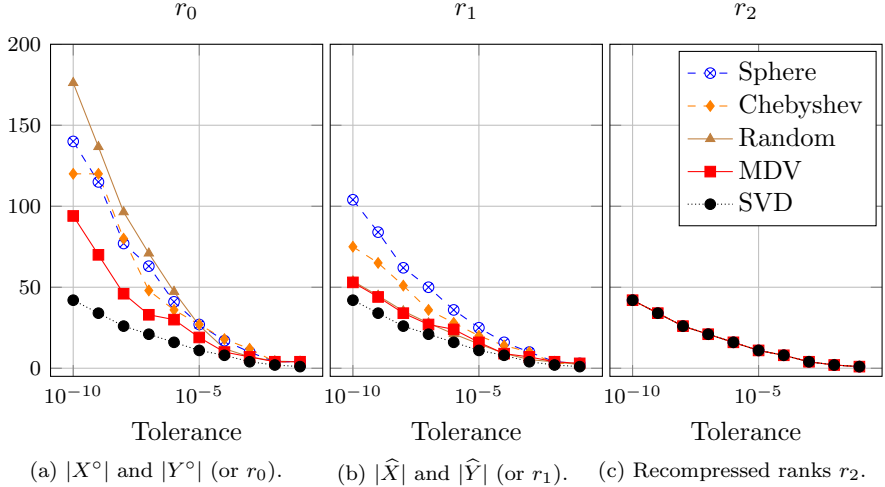
\begin{figure}[!ht]
\centering
\subfloat[$|\Xinit|$ and $|\Yinit|$ (or $r_0$). \label{fig:engine_r0}]{
\begin{tikzpicture}
\begin{semilogxaxis}[
name=plot1,
width=5.2cm,height=6cm,
ymin = -5, ymax = 200,
xlabel={Tolerance},
grid = major,
cycle list name=r0r1r2_0 list,
title={$r_0$},
]
\addplot table [x={Error}, y={r0-random}] {r0r1r2.dat};
\addplot table [x={Error}, y={r0-sphere}] {r0r1r2.dat};
\addplot table [x={Error}, y={r0-Cheb}] {r0r1r2.dat};
\addplot table [x={Error}, y={r0-mdv}] {r0r1r2.dat};
\addplot table [x={Error}, y={SVD}] {r0r1r2.dat};
\end{semilogxaxis}
\end{tikzpicture}
} \hspace{-0.5cm}
\subfloat[$|\Xhat|$ and $|\Yhat|$ (or $r_1$). \label{fig:engine_r1}]{
\begin{tikzpicture}
\begin{semilogxaxis}[
name=plot1,
width=5.2cm,height=6cm,
ymin = -5, ymax = 200,
xlabel={Tolerance},
grid = major,
yticklabel=\empty,
cycle list name=r0r1r2_1 list,
title={$r_1$},
]
\addplot table [x={Error}, y={r1-sphere}] {r0r1r2.dat};
\addplot table [x={Error}, y={r1-Cheb}] {r0r1r2.dat};
\addplot table [x={Error}, y={r1-random}] {r0r1r2.dat};
\addplot table [x={Error}, y={r1-mdv}] {r0r1r2.dat};
\addplot table [x={Error}, y={SVD}] {r0r1r2.dat};
\end{semilogxaxis}
\end{tikzpicture}
} \hspace{-0.5cm}
\subfloat[Recompressed ranks $r_2$. \label{fig:engine_r2}]{
\begin{tikzpicture}
\begin{semilogxaxis}[
name=plot1,
width=5.2cm,height=6cm,
ymin = -5, ymax = 200,
xlabel={Tolerance},
grid = major,
yticklabel=\empty,
legend entries={Sphere, Chebyshev, Random, MDV, SVD},
legend cell align={left},
cycle list name=r0r1r2_1 list,
title={$r_2$},
]
\addplot table [x={Error}, y={r2-sphere}] {r0r1r2.dat};
\addplot table [x={Error}, y={r2-Cheb}] {r0r1r2.dat};
\addplot table [x={Error}, y={r2-random}] {r0r1r2.dat};
\addplot table [x={Error}, y={r2-mdv}] {r0r1r2.dat};
\addplot table [x={Error}, y={SVD}] {r0r1r2.dat};
\end{semilogxaxis}
\end{tikzpicture}
}
\caption{Recompressing one pair of the coil case with $\dr = 1.39$}
\label{fig:r0r1r2}
\end{figure}

To illustrate the effect of the recompression step discussed in
\autoref{subsec:Discussion} on one specific pair of the coil geometry. 
\Autoref{fig:r0r1r2} shows, 
as expected, that the newly obtained low-rank approximations all have
ranks $r_2 \leq r_1$ very close to the SVD rank of $\KXY$.
Experiments on other pairs of clusters produce similar results.

\pgfplotscreateplotcyclelist{r0r1r2_0 list}{
brown, every mark/.append style={solid, fill=brown}, mark=triangle*\\
dashed, blue, every mark/.append style={solid, fill=white}, mark=otimes*\\
dashed, orange, every mark/.append style={solid, fill=orange},mark=diamond*\\%
red, every mark/.append style={solid, fill=red}, mark=square*\\%
densely dotted, black, every mark/.append style={solid, fill=black}, mark=*\\%
}
\pgfplotscreateplotcyclelist{r0r1r2_1 list}{
dashed, blue, every mark/.append style={solid, fill=white}, mark=otimes*\\
dashed, orange, every mark/.append style={solid, fill=orange},mark=diamond*\\%
brown, every mark/.append style={solid, fill=brown}, mark=triangle*\\
red, every mark/.append style={solid, fill=red}, mark=square*\\%
densely dotted, black, every mark/.append style={solid, fill=black}, mark=*\\%
}

\section{Conclusion} \label{Conclusion}

We introduced four heuristics for selecting an initial set of points for the Skeletonized Interpolation method: Chebyshev grids, Maximally Dispersed Vertices, points on a sphere, and random sampling. Some of these methods use endo-points, i.e., a subset of the given points (MDV and random), while others use exo-points, i.e., they introduce new points for the interpolation (Chebyshev and sphere).

These methods should be considered as a way to build an initial low-rank approximation at the smallest possible cost. Once a low-rank factorization exists, it can always be further compressed to have a near-optimal rank.

\textbf{SI-Chebyshev} is robust with guaranteed accuracy, even for complicated geometries and very large clusters. But it can be inefficient when points nearly lie on a submanifold.

\textbf{SI-MDV,} as a heuristic, is efficient and accurate, and very simple to implement. It may become inefficient for specific kernels, like Green's functions of Laplace equation. In that case, a large set of initial points ($r_0$) may be required.

\textbf{SI-Sphere} constructs initial points on a 2D surface instead
of a 3D geometry. This largely reduces the number of initial points,
   asymptotically. However, this method only works for kernels that
   are Green's functions of Laplace equation, and tends to
   overestimate the final rank ($r_1$). When points are distributed
   uniformly in the volume (rather than on a submanifold), {\sf
      SI-sphere} does very well for those Green's functions.

\textbf{SI-Random} is robust and general and is the easiest to implement. Nevertheless, randomly sampling the points can lead to redundancy thus making the sizes of initial points too large.

Three benchmark tests are performed, in which the four SI methods are applied to torus, plate-coil and
engine geometries. The main conclusions are summarized below:
\begin{enumerate}
  \item A comparison between \textbf{SI-Sphere and SI-Chebyshev}
  shows that the final rank approximated by {\sf SI-Chebyshev} is lower
  than that of {\sf SI-sphere} in all cases, thus numerically {\sf
     SI-Chebyshev}
is more accurate than {\sf SI-sphere}.
\item Given a pair of clusters, \textbf{SI-MDV} builds an MDV set
whose size is much smaller than the randomly sampled set constructed
by {\sf SI-random}. Therefore, {\sf SI-MDV} is more efficient than
{\sf SI-random}.
\item Compared with \textbf{SI-MDV, SI-Chebyshev} always constructs a larger size for both $\Xinit$ and $\Xhat$, thus is less
efficient and accurate than {\sf SI-MDV}. This result
suggests that in dealing with small clusters or complicated
geometries, using {\sf SI-Chebyshev} is not advantageous, despite its theoretically
guaranteed accuracy and robustness. As another special case of SI
method, {\sf SI-MDV} method is a good complement for {\sf
   SI-Chebyshev}.
\item When points are distributed uniformly inside the volume and
the kernel is a \textbf{Green's function, {\sf SI-sphere}} does very well
and is superior to {\sf SI-MDV}. The reason is that the optimal
choice of points consists in choosing points near the boundary of
the domain, while {\sf SI-MDV} samples points uniformly in the volume. This leads to overestimating $r_0$ with MDV.
\end{enumerate}

Overall, for small clusters, we recommend using \textbf{SI-MDV}, for specific Green's functions, \textbf{SI-Sphere}, for large clusters, \textbf{SI-Chebyshev}, and \textbf{SI-Random} as the simplest algorithm to implement.

\bibliography{reference}{}
\bibliographystyle{acm}

\end{document}